\documentclass[12pt]{amsart}
 \usepackage{graphicx}
\usepackage{xstring}
\usepackage{forloop}
\usepackage{bbm, amsmath, amsfonts, amscd, latexsym, amsthm, amssymb, graphicx, mathrsfs, comment}
\usepackage[colorlinks=true]{hyperref}

\usepackage{abstract}
\usepackage{mathtools}
\usepackage{caption}
\usepackage{MnSymbol}
\usepackage[dvipsnames]{xcolor}

\usepackage{tikz-cd}
\tikzset{
  symbol/.style={
    draw=none,
    every to/.append style={
      edge node={node [sloped, allow upside down, auto=false]{$#1$}}}
  }
}

\makeatletter
\newif\if@check@engine  \@check@enginetrue 
\makeatother
\usepackage[usenames,dvipsnames]{pstricks}

\usepackage{floatrowbytocbasic}

\newcommand{\nocontentsline}[3]{}
\newcommand{\tocless}[2]{\bgroup\let\addcontentsline=\nocontentsline#1{#2}\egroup}

\usepackage{booktabs}

\newtheorem{theor}{\hspace{1cm}{\sc Theorem}}[section]
\newtheorem{utver}[theor]{\hspace{1cm}{\sc Proposition}}
\newtheorem{sledst}[theor]{\hspace{1cm}{\sc Corollary}}
\newtheorem{lemma}[theor]{\hspace{1cm}{\sc Lemma}}

\newtheorem{assum}[theor]{\hspace{1cm}{\sc Assumption}}

\newtheorem*{utver*}{\hspace{1cm}{\sc Proposition}}
\theoremstyle{definition}

\newtheorem{defin}[theor]{\hspace{1cm}{\sc Definition}}
\newtheorem*{defin*}{\hspace{1cm}{\sc Definition}}
\newtheorem{exa}[theor]{\hspace{1cm}{\sc Example}}
\newtheorem{observ}[theor]{\hspace{1cm}{\sc Observation}}
\newtheorem{rem}[theor]{\hspace{1cm}{\sc Remark}}

\newtheorem{conven}[theor]{\hspace{1cm}{\sc Convention}}

\newcommand{\codim}{\mathop{\rm codim}\nolimits}

\newcommand{\sing}{\mathop{\rm sing}\nolimits}

\newcommand{\rk}{\mathop{\rm rk}\nolimits}

\newcommand{\conv}{\mathop{\rm conv}\nolimits}

\newcommand{\Trop}{\mathop{\rm Trop}\nolimits}

\newcommand{\dist}{\mathop{\rm dist}\nolimits}

\newcounter{idx}

\newcommand{\rotraise}[1]{
  \StrLen{#1}[\slen]
  \forloop[-1]{idx}{\slen}{\value{idx}>0}{
    \StrChar{#1}{\value{idx}}[\crtLetter]
    \IfSubStr{tlQWERTZUIOPLKJHGFDSAYXCVBNM}{\crtLetter}
      {\raisebox{\depth}{\rotatebox{180}{\crtLetter}}}
      {\raisebox{1ex}{\rotatebox{180}{\crtLetter}}}}
}

\renewcommand{\emph}[1]{{\it {\color{NavyBlue} #1}}}

\def\R{\mathbb R}

\def\Z{\mathbb Z}

\def\C{\mathbb C}
\def\CC{({\mathbb C}^\star)}

\def\CP{\mathbb C\mathbb P}

\def\Pr{\mathbb P}
\def\K{\mathbb K}

\usepackage{xr-hyper}
\makeatletter

\newcommand*{\addFileDependency}[1]{\typeout{(#1)}

\@addtofilelist{#1}
\IfFileExists{#1}{}{\typeout{No file #1.}}
}\makeatother

\begin{document}

\begin{center}{\Large \sc Whitney fold\&cusp for algebraic surfaces, and singularities of discriminants}

\vspace{3ex}

{\sc Alexander Esterov\footnote{LIMS, aes@lims.ac.uk}, Lev Vladimirov\footnote{HSE, Vladimirov.L.S@hse.ru}, Aliaksandr Yuran\footnote{HSE, ayuran@hse.ru}}
\end{center}

\begin{abstract}
We describe transversal singularity types of the singular locus of the $A$-discriminant for $A\subset\Z$, and deduce a Whitney type theorem for a coordinate projection of a surface defined by a general polynomial equation with a given Newton polytope $N$: under mild combinaorial conditions on $N$, all multisingularities are stable (folds, cusps, and double folds). We then enumerate the multisingularities in terms of $N$. The results rely on the analysis of degeneracy of relevant Vandermonde/Schur type matrices, which may be of independent interest.
\end{abstract}

\tableofcontents

\section{introduction}

{\bf The question: stable singularities of regular maps.} A {\it singularity} is an equivalence class of germs of smooth maps $(P,p)\to(Q,q)$ under a certain type of local authomorphisms of the manifolds $(P,p)$ and $(Q,q)$ -- most notably, self-diffeomorphisms or self-homeomorphisns. The latter choice does not matter over $\C$, when it comes to simplest singularities, such as
$$(\K^2,0)\to(\K^2,0),\quad (x,y)\mapsto (x,y^2) \quad \mbox{({\it fold}) and}$$
$$(\K^2,0)\to(\K^2,0),\quad (x,y)\mapsto (x,y^3+xy)  \quad \mbox{({\it cusp}).}$$
These two are notable for being {\it stable}: perturbing a representative of such a singularity gives a map with the same singularity at a nearby point.

A {\it multisingularity} is an equivalence class of finite tuples of maps $(P_i,p_i)\to(Q,q)$, sharing the common range. It can be stable as well:
$$(\K^2,0)\sqcup(\K^2,0)\to(\K^2,0),\quad (x_1,y_1)\mapsto (x_1,y_1^2),\quad (x_2,y_2)\mapsto (x_2^2,y_2)  \quad \mbox{({\it double fold}).}$$

\begin{figure}[hbt!]
\floatbox[\capbeside\thisfloatsetup{capbesideposition={right,center}}]{figure}[7cm]
  {
    \caption{The line is the discriminant of a projection of a torus into $\R^2$. The blue and red dots are the cusps and double folds, the rest are single folds.}
  }  {\includegraphics[width=5cm]{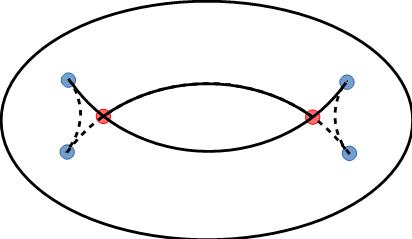}
  }
\end{figure}

In \cite{??whitney}, Whitney proved there are no other stable (multi)singularities in this dimension, and any smooth proper map of real surfaces can be perturbed to a {\it stable} one, i.e. having only cusps, folds and double folds (modulo removing regular germs from multigerms).
Is it true for general maps of algebraic surfaces?

The answer depends on what family of maps we take a general element from. For instance, if the map $\C^2\to\C^2$ is given by a pair of general polynomials with prescribed Newton polygons, 
the answer depends on the polygons and is difficult (unless the polygons have a large common part: e.g. see \cite{jel} for maps given by general polynomials of prescribed degrees, and \cite{hilani} for how much less we know for arbitrary polygons).

To gain understanding, we study this question in a simpler setting of the same kind (see the next remark 
for why we deem it simpler): given the projection $p$ of a general surface $V\subset\CC^3$ with a prescribed Newton polytope $N\subset\R^3$ to the first coordinate plane $\CC^2$, what singularities could it have, and how many?

One natural approach is to extend the domain of the projection $p$ to the $N$-toric compactification of the surface $V$, and count sigularities of the resulting proper map using Thom polynomials. This works well if $N=d\cdot($standard simplex$)$, so that $V$ is a general surface of degree $d$, and the $N$-toric compactification of $\CC^3$ is $\CP^3$, see \cite{ohmoto} for details.
But this fails tor most of polytopes $N$, because the proper extension of the projection acquires complicated non-stable singularities at the toric boundary.

\vspace{1ex}

{\bf Reduction to singularities of sparse discriminants.} Gelfand--Kapranov--Zelevin\-sky's {\it $A$-discriminant} $D_A$ is the closure of all univariate polynomials with a multiple non-zero root in the space $\C^A:=\{\sum_{a\in A}c_at^a\},\,A\subset\Z$ \cite{gkz}. It is obviously related to the projection $p$ of our interest for $A:=h(N\cap\Z^3)$, where the {\it height} $h:\Z^3\to\Z$ denotes the third coordinate. Indeed, the map $\CC^2\to\C^A$, sending each point $(x_1,x_2)$ to the fiber $f(x_1,x_2,\cdot)$, sends the discriminant of $p$ to $D_A$. This translates the initial question to a version for $A$-discriminants: what  singularity does $D_A$ have at a general point of its singular locus? 

\begin{rem}
The same trick applies to a map $\CC^2\to\C^2$ given by general polynomials with Newton polygons $P$ and $Q$. This will lead to the $(P,Q)$-discriminant of two bivariate polynomials, drastically more complicated than $D_A,\, A\subset\Z$.
\end{rem}

In the case $A=\{0,1,\ldots,d\}$, the hypersurface $D_A$ is the classical Sylvester discriminant. Its singular locus is the closure of two strata:

$\prec\; :=\{$polynomials with a triple root$\}$, with a cusp as transversal singularity, and 

$\times:=\{$polynomials with two double roots$\}$, with a node as transversal singularity.

(The {\it transversal singularity} of a component $S\subset\sing D_A$ is the singularity of $D_A\cap($a general germ of a 2-plane intersecting $S)$, if it is the same for almost all such germs.)

To what extent is it true for arbitrary $A\subset\Z$?

\begin{conven}\label{conven1} We always assume the support set $A$ has at least six elements, and no proper arithmetic progression $a\cdot\Z+b\subsetneq\Z$ contains $|A|-1$ of its elements.\end{conven}

Our results and proofs adapt to support sets breaking the progression condition: we do not cover this in the paper because the adaptation is straightforward but tiresome. Support sets with less than six points, on the contrary, demonstrate genuinely sporadic behavior requiring  separate analysis.
\begin{theor}[\ref{mainth}]\label{thsingd} Let the support set $A\subset\Z$ satisfy Convention.
The $A$-discriminant $D_A$ has only cusps and nodes as transversal singularities at general points of its singular locus components of codimension 2, if \textbf{\textit{and only if}} $A\subset\Z$ contains $\{1+\min A,3+\min A\}$ or $1+\min (A\setminus\min A)$, and contains $\{\max A-1,\max A-3\}$ or $\max (A\setminus\max A)-1$. 

\end{theor}
For the proof and exact description of the singularity strata, see Section 3.
\begin{rem}\label{rem1}
1. The singular locus of the discriminant of $f(t)=a_0+a_1t+a_4t^4+a_8t^8$ has a component $a_0=a_1=0$ whose transversal singularity is $(s^3,s^4)$, and a component $a_1=a_4^2-4a_0a_8=0$ whose transversal singularity is a union of four lines. 

2. This theorem may be strengthened in several directions beyond the scope of this work. There may be at most three cuspidal and arbitrarily many nodal components of $\sing D_A$ (Remark \ref{wildcomponents}). Dropping the assumption on nearly extremal points on $A$, some of the components of $\sing D_A$ may exhibit more complicated irreducible singularities: their Milnor numbers are computed in Section 3, and the method actually allows to compute all characteristic Puiseaux exponents (determining the topological type of the singularity). One can show there are no components of codimension higher than 2 in $\sing D_A$, though it is not necessary for our application to singularities of general surface projections.

3. There are other motivations to study generic singularities of sparse resultants and discriminants, see e.g. \cite{weze, mts}. That they are not actively studied is probably because of seemingly hopeless difficulty of such questions, rather than lack of motivation. We hope that the surprisingly complete answer in Theorem \ref{thsingd} is reassuring in this respect.

4. Singularities of sparse discriminants and resultants may look difficult to study, because the standard count of infinitesimal deformations reduces the question to ranks of Vandermonde/Schur type matrices like $M=(x_i^{a}),\,a\in A\subset\Z$, which are notoriously complicated \cite{shapiro}. We shall see however that the questions about such matrices coming from singularities of discriminants sometimes admit surprisingly full and interesting answers.
The ones that we need are derived in Section 2.
\end{rem}
\begin{theor}[\ref{mainrank}]
There are finitely many points $(p_1:p_2:p_3)\in\CC^3/\CC, p_i\ne p_j$, such that the $3\times |A|$ matrix $(p_i^a)_{a\in A,\; i=1,2,3}$ is degenerate. For every such point,  $$\rk (p_1^a,p_2^a,p_3^a,ap_1^a,ap_2^a,ap_3^a)_{a\in A}\ge 4.$$ 
\end{theor}

{\bf Singularities of projections.} 
We call $N\subset\Z^3$ a {\it polytope} if it is the intersection of a convex polytope in $\R^3$ with $\Z^3$, and call its subset a {\it face}, if it is its  intersection with a face of its convex hull. 
Let $f:\CC^3\to\C$ be a general Laurent polynomial with the Newton polytope $N$. We study the projection $p:\{f=0\}\to\CC^2$ forgetting the third coordinate, and call the third coordinate $h:\Z^3\to\Z$ the {\it height}. We assume $|h(N)|>2$.

\begin{theor}\label{thsingp} If $h(N)\subset\Z$ satisfies Convention \ref{conven1}, then all multisingularities of the surface projection $p$ are stable (folds, cusps and double folds, modulo removing regular germs from multigerms).
\end{theor}
Note that the statement about monosingularities being only folds and cusps requires a much weaker assumption $|h(N)|>2$.

To express the number of singularities, 
denote the convex hull by $[\cdot]$, and define the {\it incremental polytopes} of $N$ as 

$$N':=\bigcap_{q\in\Z} [N+(N\setminus h^{-1}q)]=\bigcap_{q\in\Z} [(N+N)\setminus h^{-1}q]=[a+b\,|\,a,b\in N,\,h(a)\ne h(b)];\eqno{(*')}$$

$$N'':=\bigcap_{q_1\ne q_2} [N+(N\setminus h^{-1}q_1)+(N\setminus h^{-1}\{q_1,q_2\})]=\bigcap_{q} [N+(N\setminus h^{-1}q)']=$$
$$=[a_1+a_2+a_3\,|\,a_i\in N,\,h(a_i)\mbox{ are pairwise different}].\eqno{(*'')}$$

\begin{figure}[hbt!]
\includegraphics[width=10cm]{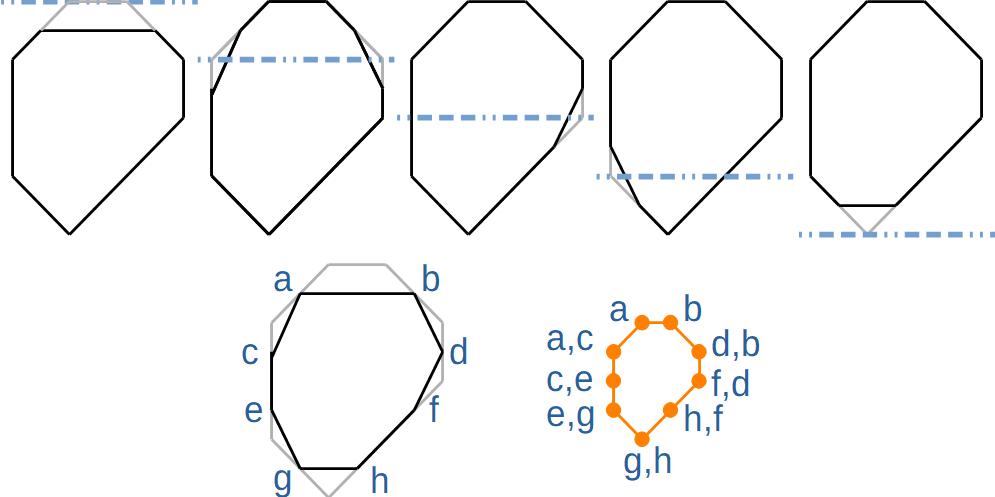}

\caption{Computing the incremental polytope $N'$ for the orange $N$: \newline
the top row of polygons is what intersects in the 2nd expression of $(*')$;\newline
representing each vertex of $N'$ as a sum of two points of $N$ at different height (denoted by the same letter) illustrates the third expression of $(*')$.}
\end{figure}

\begin{theor}[\ref{thtbenum} for $\dim=3$]\label{enumcusp}
1. The projection $p$ has $N\cdot(N'-N)\cdot(N''-N')$ cusps (to evaluate this expression, open the brackets, and evaluate every monomial $P\cdot Q\cdot R$ as the lattice mixed volume of the polytopes $P,Q$ and $R$).

2. The Euler characteristic of the critical locus of $p$ (which is the set of all cusp and fold points) equals $N\cdot N'\cdot(N-N')$.

3. The tropical fan of the critical locus equals $N\cdot(N'-N)=N\cdot N'-N\cdot N$ (evaluating the product of polytopes as the intersection product of their dual fans).
\end{theor}

Geometrically, the $k$-th incremental polytope $N^{(k)}$ is a polytope whose support function is the sum of tropicalizations of the logarithmic derivatives $f^{(i)}|_{f^{(i-1)}=\cdots=f=0}$ over $i=0,\ldots,k$; see \cite{eci} or Section 4 for detail. With this geometric definition of $N^{(k)}$, the theorem is a special case of the BKK formula for sch\"on complete intersections \cite{sci} (they are a common generalization for toric complete intersections, hyperplane arrangement complements, and many other special algebraic varieties of combinatorial nature).
The first two expressions for the incremental polyope $N^{(k)}$ in $(*')$ and $(*'')$ were found in \cite{eci}, and the thrid one in \cite{kks}. 

\begin{rem} 1. The first and the third expression for the incremental polytopes involve certain Minkowski pointwise summation, which makes them difficult to compute and estimate in practice. The second expression is free from this flaw, but is much more difficult to establish. 

2. The classical Bernstein--Kouchnirenko--Khovanskii formula \cite{Be, Kouch, Kh}, expressing the geometry of a toric complete intersection in terms of the Newton polytopes of its equations, is not applicable to the sch\"on complete intersection of our interest $f^{(k)}=\cdots=f'=f=0$ even for general $f$: that is why we need its generalization \cite{sci}.

3. For most of polytopes $N$, the critical locus of $p$ is an irreducible curve (\cite{sci} and \cite{zh} give some specific criteria for irreducibility of sch\"on complete intersections). In this case, the genus of the critical curve can be extracted from Theorem \ref{enumcusp} as one half of $2-($the total multiplicity of the rays of the tropical fan$)-($the Euler characteristics$)$.

4. For most polytopes $N$, the discriminant curve of the projection $p$ has only nodes and cusps as singularities. For example, if $N$ is the standard simplex of size 3 (i.e. $f=0$ is a general cubic surface), the discriminant curve of its coordinate projection has degree 6 with only simple cusps and nodes. Classifying possible topological types of the real part of this curve for real $f$ is an open question in the vein of Hilbert's 16'th problem. 

Another natural open question is whether, for any degree $d$, there exists a real algebraic surface of degree $d$, all of whose projection's cusps are real. The positive answer is classically known for $d=3$ and is computer-verified for $d=4$ in \cite{emm}.
\end{rem}

In Section 4, we prove a multidimensional version of Theorems \ref{thsingp} and \ref{enumcusp} for the Thom--Boardman singularity strata of a coordinate projection of a general hypersurface in $\CC^n$, and give several ways to count the double folds of $p$ in terms of the Newton polytope $N$, similarly to the above formula for cusps. This requires computing a so call ultratropicalization of the discriminant curve, which is done in Section 5.

Sections 2 and 3 are devoted to ranks of Vandermonde type matrices, and to the consequences for singularities of discriminants, respectively.

\vspace{1ex}

A significant part of Sections 2 and 3 is based on the BSc diploma work \cite{lv}. ChatGPT was used for bookkeeping and verification of computationally intensive proofs at the late stages of the work.

\section{Ranks of Vandermonde type matrices}
\begin{defin}
    Consider a set $A=\{a_1, \dots, a_n\}\subset \mathbb Z$. The set $A$ is reduced if it cannot be shifted into a proper sublattice $k\mathbb Z \subset \mathbb Z$. 
    
    By $\gcd (A-A)$ denote the greatest common divisor of all differences of the form $a_i-a_j$. Then $A$ is reduced if and only if $\gcd(A-A)=1$.

    We call he set $A$ doubly reduced if the set $A\setminus a_i$ is reduced for each $i$.
\end{defin}

\begin{defin}
	 	Consider a set $A=\{a_1, \dots, a_n\}\subset \mathbb Z$. Define the $m \times n$ matrix for any nonzero complex $p$ and positive integer $m$:
	\begin{equation*}
		M_{A}^m(p)=\begin{bmatrix}
			p^{a_1} & p^{a_2} &  & p^{a_n} \\
			a_1 p^{a_1} & a_2 p^{a_2} & \dots & a_n p^{a_n} \\
			a_1^2 p^{a_1} & a_2^2 p^{a_2} &  & a_n^2 p^{a_n} \\
			\vdots & \vdots & \vdots & \vdots \\
			a_1^{m-1} p^{a_1} & a_2^{m-1} p^{a_2} &  & a_n^{m-1} p^{a_n}
		\end{bmatrix}.
	\end{equation*}

Let $S=(m_1,\dots,m_k)$ be a sequence of positive integers, and $\mathbf{p}=(p_1,\dots,p_k)$ be a sequence of distinct points in $\C^*$. Define the generalized Vandermonde matrix $M_{A}^S(p_1,\dots,p_k)$ of size $\left(\sum m_i\right) \times n$ by concatenation of matrices:

	\begin{equation*}
		M_{A}^S(p_1,\dots,p_k)=\begin{bmatrix}
			M_{A}^{m_1}(p_1) \\
			M_{A}^{m_2}(p_2) \\
			\vdots  \\
			M_{A}^{m_k}(p_k)
		\end{bmatrix}.
	\end{equation*}
\end{defin}

\begin{exa}
	The matrix $M_{\{0,1,2,\dots,k-1\}}^{1,1,\dots,1}(p_1,\dots,p_k)$ is the classical Vandermonde matrix. If $p_1,\dots,p_k \in \mathbb Z$, then $M_{\{0,1,2,\dots,k-1\}}^{1,1,\dots,1}(p_1,\dots,p_k)=M_{\{p_1,\dots,p_k\}}^{k}(1)^T$.

    The fraction $s_{a_1,\dots,a_k} = \det M_{a_1,\dots,a_k}^{1,\dots,1}(p_1,\dots,p_k)/\det M_{\{0,1,\dots,k-1\}}^{1,\dots,1}(p_1,\dots,p_k)$ is a classical Schur polynomial (note that the indices differ from the standard notation).
\end{exa}

\begin{observ} \label{vandermonderank}
    Consider a polynomial $f = \sum_{a\in A} c_a t^a$ with roots at points $p_1 \dots p_k \in \C^* $ of multiplicities $m_i+1$ respectively and no other roots of multiplicity greater than one. 

    Consider the natural linear map $\pi \colon \C^A \to \C[t^{\pm1}]/(f,f')$. Then $\rk \pi = \rk M_{A}^{m_1,\dots,m_k}(p_1,\dots,p_k)$.
\end{observ}

\begin{utver} \label{mult}
    For any $A\subset \mathbb Z$, positive integer $m$ and $p \in \C^*$, we have $\rk M_{A}^{m}(p) = \min \{m,\lvert A \rvert \}$
\end{utver}

\begin{proof}
    If $m=\lvert A \rvert$, then $\rk M_{A}^{m}(p) = \min \{m,\lvert A \rvert \}$ is the ordinary Vandermonde matrix with columns multiplied by powers of $p$ and thus it is nondegenerate. The other cases easily follow from this one.
\end{proof}

\begin{lemma} \label{matrixtorusaction}
    Assume that $q\in \C^*$. Then \[ \rk M_A^S(p_1,\dots,p_k) = \rk M_A^S (qp_1,qp_2,\dots,qp_k).\]
\end{lemma}

\begin{proof}
    The second matrix is obtained from the first by multiplication of columns by powers of $q$. 
\end{proof}

\begin{utver} \label{oneone}
    Assume that $A$ is reduced, $\lvert A \rvert \ge 2$, and $p_1 \ne p_2 \in \C^*$. Then $\rk M_A^{1,1}(p_1,p_2) = 2$.
\end{utver}

\begin{proof}
    Assume that the two rows of the matrix $M_A^{1,1}(p_1,p_2)$ are linearly dependent. Then there exists $q\in \C$ such that for any $i$ the equality $p_1^{a_i} = qp_2^{a_i}$ holds.
    
    Since $A$ cannot be shifted to a proper sublattice, there exist integers $k_i$ such that $k_1+k_2+\dots + k_n = 0$, $k_1a_1 + k_2a_2 + \dots + k_n a_n = 1$ (indeed, there exists a $\mathbb Z$-linear combination of the numbers $a_1-a_i$ equal to 1). Thus, $p_1/p_2=\prod_{i=1}^n (p_1/p_2)^{k_i a_i} = \prod_{i=1}^n q^{k_i} = 1$.
\end{proof}

\begin{utver} \label{oneoneroots}
     Assume $\rk M_A^{1,1}(p_1,p_2) = 1$ and $p_1\ne p_2$. Then there exists an integer $k \ge 2$ such that $\gcd(A-A)=k$, and $(p_1/p_2)^k=1$.
\end{utver}

\begin{proof}
    By $k$ denote the number $\gcd(a_1-a_2,a_1-a_3,\dots,a_1-a_n)$. Then there exists a $\mathbb Z$-linear combination of numbers $a_1-a_i$ equal to $k$. Thus, there exist integers $k_i$ such that $k_1+k_2+\dots + k_n = 0$, $k_1a_1 + k_2a_2 + \dots + k_n a_n = k$. Now, $(p_1/p_2)^k=\prod_{i=1}^n (p_1/p_2)^{k_i a_i} = \prod_{i=1}^n q^{k_i} = 1$. If $p_1\ne p_2$, then $k \ge 2$.
\end{proof}

\begin{lemma}\label{onetwoabsonelemma}
    Assume that $A = \{a_1, a_2,a_3\}$ and $\lvert y \rvert=1$. Then $\rk M_{A}^{2,1}(1,y)<3$ only if $y^{a_1} = y^{a_2} = y^{a_3}$.
\end{lemma}

\begin{proof} 
    Assume that \[ \det \begin{bmatrix}
1 & 1 & 1  \\
a_1  & a_2 & a_3 \\
y^{a_1} & y^{a_2} & y^{a_3}
\end{bmatrix} 
 = (a_1-a_2)(y^{a_2}-y^{a_3})-(a_2-a_3)(y^{a_1}-y^{a_2}) = 0.\] If $y^{a_i}\neq y^{a_j}$ for $i \neq j$, then $\frac{a_1-a_2}{a_2-a_3} = \frac{y^{a_1}-y^{a_2}}{y^{a_2}-y^{a_3}}$. Since the three points $y^{a_i}$ lie on the unit circle and are pairwise different, they are not collinear. Thus, $ \frac{y^{a_1}-y^{a_2}}{y^{a_2}-y^{a_3}}$ is not a real number. However, $\frac{a_1-a_2}{a_2-a_3}$ is. Hence, without loss of generality, we may assume that  $y^{a_1} = y^{a_2}$. Now the assumption that the determinant is zero implies that $y^{a_2} = y^{a_3}$ as well.
\end{proof}

\begin{utver} \label{onetwoabsone}
    Assume that $\gcd(A-A) = 1$ and $\lvert A \rvert \ge 3$. Let $p_1,p_2 \in \C^*$ be two different numbers such that $\lvert p_1 \rvert = \lvert p_2 \rvert$. Then $\rk M_A^{1,2}(p_1,p_2) = 3$.
\end{utver}
\begin{proof}
    Assume that $\rk M_A^{1,2}(p_1,p_2) < 3$. Denote $q=p_1/p_2$. Note that $\lvert q \rvert = 1$. By Lemma \ref{matrixtorusaction}, $\rk M_A^{1,2}(1,q) < 3$. Applying the Lemma \ref{onetwoabsonelemma} to the $3\times3$ minors of the matrix $M_A^{1,2}(1,q)$, we obtain the equality $q^{a_1} = q^{a_2} = \dots = q^{a_n}$. Thus, $\rk M_A^{1,1}(1,q) = 1$. By Lemma \ref{oneone}, $q=1$.
\end{proof}

\begin{lemma} \label{fewnomialbound}
	Let $f_1,f_2\in \mathbb{R}[x]$ be two polynomials of degree $d_1$ and $d_2$ respectively and $q_1, q_2$ be two positive real numbers. Assume that $q_1\ne q_2$ or $f_1 \ne -f_2$. Then the equation $q_1^x f_1(x) + q_2^x f_2(x)=0$ has at most $d_1+d_2+1$ real solutions.
\end{lemma}

\begin{proof}
	This statement is a variant of the Descartes rule of signs. The proof is by induction in $d_1$ and $d_2$. The base case of induction is $d_1=0$ or $d_2=0$. The step $d_2-1 \to d_2$ is an application of Rolle's theorem. The function $(q_1/q_2)^x f_1(x) + f_2(x)=0$ has at most $d_1 + d_2 +1$ real solutions, since the equation $(\log(q_1/q_2) + f_1'(x)) (q_1/q_2)^x + f_2'(x)=0$ has at most $d_1 + d_2$ solutions.
\end{proof}

\begin{utver}
\label{multmult} 
	Let $p_1,p_2$ be two nonzero complex numbers such that $\lvert p_1 \rvert \ne \lvert p_2 \rvert$. Let $m_1,m_2$ be two positive integers and $A$ be a subset of $\mathbb Z$ such that $\lvert A \rvert > 2(m_1+m_2)-3$. Then \[\rk M_A^{m_1,m_2}(p_1,p_2)=m_1+m_2.\] 
\end{utver}

\begin{proof}
	Assume that $\rk M_A^{m_1,m_2}(p_1,p_2)<m_1+m_2$. Then there exists a nonzero vector \[v=(\alpha_0, \alpha_1, \dots, \alpha_{m_1-1}, \beta_0, \beta_1, \dots, \beta_{m_2-1})\in \C^{m_1+m_2}\] such that $v \cdot M_A^{m_1,m_2}(p_1,p_2)=0$. This equality is equivalent to the following condition. Each element $t \in A$ is a solution of the equation \begin{equation*}
		\alpha(t) p_1^t + \beta(t) p_2^t = 0,
	\end{equation*}  where $\alpha(t) = \sum_i \alpha_i t^i$, $\beta(t) = \sum_i \beta_i t^i$.
	
	If this equation is satisfied for each element $t\in A$, then for each $t\in A$ the equality \[\alpha(t) \overline \alpha(t) \lvert p_1^2 / p_2^2 \rvert ^t - \beta(t) \overline \beta(t) = 0.\]
	holds. By Lemma \ref{fewnomialbound} this equality holds for at most $2(m_1-1+m_2-1)+1$ real values of $t$.
\end{proof}

\begin{utver} \label{oneoneone}
    Assume that $\gcd(A-A) = 1$ and $\lvert A \rvert \ge 4$. Then the number of pairs $p_1,p_2 \in \C^*$ such that $1\ne p_1\ne p_2 \ne 1$ and $\rk M_A^{1,1,1}(1,p_1,p_2) < 3$ is finite.
\end{utver}

\begin{proof}
    This fact follows from the irreducibility of Schur polynomials \cite{dz}.
    
    Consider the variety $X=\{(p_1,p_2) \in (\C^* \setminus 1)^2 \mid p_1 \ne p_2, \rk M_A^{1,1,1}(1,p_1,p_2) < 3\}$. Our aim is to show that $X$ is finite.
    
    Consider the variety $X' \subset \Pr^2$ obtained as the intersection of all curves $X'_{ijk}$ $(a_i<a_j<a_k)$ defined by equations $s_{a_i,a_j,a_k}(p_1,p_2,p_3)=\frac{\det M_{a_i,a_j,a_k} (p_1,p_2,p_3)}{  (p_1-p_2)(p_1-p_3)(p_2-p_3)} = 0$ for $\{a_i,a_j,a_k\} \subset A$. Then $X$ is the intersection of $X'$ and the affine chart $p_3=1$.

    Note that $s_{a_i,a_j,a_k}(p_1,p_2,p_3)$ is a Schur polynomial. Denote $v_{a_i,a_j,a_k}(p_1,p_2,p_3)=\frac{(p_1^d-p_2^d)(p_1^d-p_3^d)(p_2^d-p_3^d)}{(p_1-p_2)(p_1-p_3)(p_2-p_3)}$, here $d = \gcd(a_i-a_j,a_i-a_k)$. Then by \cite{dz} the Schur polynomials $s_{a_i,a_j,a_k}$ can be represented in the form $s_{a_i,a_j,a_k} =(p_1 p_2 p_3)^{-a_i} \tilde s_{a_i,a_j,a_k} v_{a_i,a_j,a_k}$ where $\tilde s_{a_i,a_j,a_k}$ is an irreducible polynomial. 

    It is enough to show that the greatest common divisor of the polynomials $s_{a_i,a_j,a_k}$ in the ring $\C[p_1,p_2,p_3]$ is $1$, that is \begin{equation*}
        \gcd_{\{a_i,a_j,a_k\}\subset A} \left(\tilde s_{a_i,a_j,a_k}(p_1,p_2,p_3) v_{a_i,a_j,a_k}(p_1,p_2,p_3)\right)
    \end{equation*} First, note that the polynomials $v_{a_i,a_j,a_k}$ have no common factors. Indeed, assume that $p_I-\varepsilon p_J$ is their common divisor ($\varepsilon$ is a $D$-th root of unity, $I,J\in \{1,2,3\}$). Then for any $i,j$ the difference $a_i-a_j$ is divisible by $D$, which contradicts the condition $\gcd(A-A)=1$. Moreover, $p_I-\varepsilon p_J$ does not coincide with any of the irreducible symmetric polynomials $\tilde s_{a_i,a_j,a_k}$. Hence, \begin{equation*}
        \gcd_{\{a_i,a_j,a_k\}\subset A} \left(\tilde s_{a_i,a_j,a_k}v_{a_i,a_j,a_k}\right)=\gcd_{\{a_i,a_j,a_k\}\subset A} \left(\tilde s_{a_i,a_j,a_k}\right)
    \end{equation*}

    Now, it remains to show that some two of the polynomials $\tilde s_{a_i,a_j,a_k}$ do not coincide. Indeed, let $A'=\{a,a+p,a+p+q,a+p+q+r\}$ be a subset of $A$ such that $\gcd(A'-A')$ is not divisible by 3.
    
      Note that $\deg \tilde s_{a_i,a_j,a_k} = a_k+a_j-2a_i-3\gcd(a_j-a_k,a_i-a_j).$ Consider the degrees of the 4 polynomials $\tilde s_{a_i,a_j,a_k}$ corresponding to the subsets of three elements of $A'$. They are equal to $2p+q$, $2p+q+r$, $2p+2q+r$, $2q+r$ modulo 3. If all these remainders coincide, then the numbers $p,q,r$ are divisible by 3. But $\gcd(A'-A')$ is not divisible by 3, so two of the polynomials $\tilde s_{a_i,a_j,a_k}$ have different degrees and therefore do not coincide.
\end{proof}

\begin{utver}[a small generalization of \cite{esv}, Lemma 5.22] \label{oneoneoneabsone}
    Let $A$ be a support such that $\gcd(A-A) = 1$ and $p_1,p_2$ be different numbers not equal to $1$. Assume that $\lvert p_1 \rvert = \lvert p_2 \rvert = 1$ and $\rk M_{A}^{1,1,1}(1,p_1,p_2)<3$. Then there exists an integer $k\ge 3$ such that $p_1^k=p_2^k=1$. Moreover, $A$ splits into two nonempty subsets $A = B_1 \sqcup B_2$  such that the numbers $\gcd(B_1-B_1)$ and $\gcd(B_2-B_2)$ are divisible by $k$.  
\end{utver}

\begin{proof}
    Since $M_{A}^{1,1,1}(1,p_1,p_2)$ is degenerate, consider a nonzero vector $(S,T,U)$ such that $(S,T,U) \cdot M_A^{1,1,1} = 0$. Thus, for any $a_i\in A$ the equality $S+Tp_1^{a_i}+Up_2^{a_i}=0$ holds.
    
    If $S=0$, then $\rk M_A^{1,1}(p_1,q_1) \le 1$ which contradicts Lemme \ref{oneone}. Thus, $S\ne 0$. Similarly, $T\ne 0$ and $U \ne 0$. Consider the affine transformation $V \colon \C \to \C$, $x \mapsto -S/U - (T/U)x$. Then $V(p_1^{a_i})=p_2^{a_i}$ for all $a_i \in A$.

    Consider the set $P=\{p_2^{a} \mid a\in A \}$. Assume that $\lvert P \rvert \ge 3$. Then choose three distinct elements $p_2^{a_i}$, $p_2^{a_j}$, $p_2^{a_k}$. Then the transformation $V$ sends 3 points $p_1^{a_i}$, $p_1^{a_j}$, $p_1^{a_k}$  on the unit circle to 3 distinct points on the unit circle. Hence, the center of the circle is mapped to itself: $V(0)=0$ which means that $S=0$. But $S\ne 0$, so $\lvert P \rvert \le 2$.

    Hence, $A$ is the union of two subsets $A = B_1 \sqcup B_2$ such that $p_2^{a_i} = p_2^{a_j}$ if $a_i$ and $a_j$ belong to the same subset. This means that $\rk M_{B_i}^{1,1}(1,p_2) = 1$ for $i=1,2$. Proposition \ref{oneoneroots} implies that there exists $k\ge 2$ such that $p_2^k=1$, $\gcd(B_1-B_1)$ and $\gcd(B_2-B_2)$ are divisible by $k$. Similarly, we obtain $p_1^k=1$. If $k=2$, then the condition $1\ne p_1\ne p_2 \ne 1$ can not be satisfied.
\end{proof}

\begin{lemma} \label{mobius}
    Assume that $\rk M_A^{2,2} (1,p) < 4$. Then one of the following statements holds.
    \begin{itemize}
        \item There exists a nonconstant M\"obius transformation sending the points $a_1,\dots,a_n$ to $p^{a_1},\dots,p^{a_n}$ respectively. In particular, the points $p^{a_1},\dots,p^{a_n}$ are distinct.
        \item Some $n-1$ of the points $p^{a_1},\dots,p^{a_n}$ coincide. 
    \end{itemize}
\end{lemma}
\begin{proof}
    There is a nonzero vector $(S,T,U,V)$ such that $(S,T,U,V)  \cdot M_A^{2,2} (1,p) = 0$, that is, for any $i$ the equality $S+Ta_i + Up^{a_i} + Vp^{a_i}a_i=0$ holds. 
    
    Assume that for any $i$, the number $U+Va_i$ is nonzero. Then the map $x \mapsto \frac{-S-Tx}{U+Vx}$ is the desired one. Note that this map may be constant if $SV=TU$. Then $p^{a_1}=p^{a_2}=\dots=p^{a_n}$.

    Now, assume that $U+Va_j=0$ for some $j$. Then $S+Ta_j=0$.
    
    Note that the first two rows of $M_A^{2,2} (1,p)$ are linearly independent, so at least one of the numbers $U$ and $V$ is nonzero. It now follows that $S=\lambda U$, $T = \lambda V$ for some $\lambda \in \C$. Indeed, the system of equations $U+xV=S+xT=0$ has the unique solution $x=a_j$. Note that $U+Va_i \ne 0$ for $i\ne j$.

    Hence, $p^{a_i} = \frac{-S-Ta_i}{U+Va_i} = -\lambda$ for any $i\ne j$.
\end{proof}

\begin{lemma} \label{twotwounity}
    Assume that $A$ is doubly reduced. Let $a_j$ and $a_k$ be any two distinct elements of $A$ and $p\ne 1$ be a $(a_k-a_j)$-th root of unity. Then $\rk M_A^{2,2}(1,p) = 4$.
\end{lemma}

\begin{proof}
    Assume that $\rk M_A^{2,2}(1,p) < 4$. Note that $p^{a_k} = p^{a_j}$. By Lemma \ref{mobius}, there exists $a_i \in A$ such that the value $p^{a}$ does not depend on the choice of an element $a\in A \setminus \{a_i\}$. 

    Thus, $\rk M_{A \setminus \{a_i\}}^{1,1}(1,p)=1$. By Lemma \ref{oneone}, the set $A \setminus \{a_i\}$ can be shifted into a proper sublattice of $\mathbb Z$.
\end{proof}

\begin{theor} \label{mainrank}
    Assume that $\lvert A \rvert \ge 6$ and that $A$ is doubly reduced. Then the following statements hold.
    
    (1) If $m \in \mathbb Z$, $1\le m \le \lvert A \rvert$ and $p\in \C^*$, then $\rk M_A^{m}(p) = m$.

    (11) If $p_1,p_2 \in \C^*$, $p_1\ne p_2$, then $\rk M_A^{1,1} (p_1,p_2) = 2$.

    (21) If $p_1,p_2 \in \C^*$, $p_1\ne p_2$, then $\rk M_A^{2,1} (p_1,p_2) = 3$.

    (111) The set $X_A=\{(p_1,p_2)\in \CC^2 \mid \rk M^{1,1,1}(1,p_1,p_2) < 3, 1\ne p_1 \ne p_2 \ne 1 \}$ is finite.

    (222) If $(p_1,p_2)\in X_A$, then $\rk M^{2,2,2}(1,p_1,p_2) \ge 4$.
\end{theor}

\begin{proof}
    (1) is Lemma \ref{mult}.
    
    (11) is Lemma \ref{oneone}

    (21) is Lemma \ref{multmult} if $\lvert p_1 \rvert \ne \lvert p_2 \rvert$ and Lemma \ref{onetwoabsone} if $\lvert p_1 \rvert = \lvert p_2 \rvert$.

    (111) is Lemma \ref{oneoneone}

    (222) Assume that $\lvert p_i \rvert \ne 1$ for $i=1$ or $i=2$. Then by Lemma \ref{multmult} $\rk M^{2,2}(1,p_i) = 4$, so $\rk M^{2,2,2}(1,p_1,p_2)$ contains a submatrix of rank $4$.

    Now, assume that $\lvert p_1 \rvert = \lvert p_2 \rvert = 1$. Since $(p_1,p_2)\in X_A$, Lemma \ref{oneoneoneabsone} implies that $A=B_1 \sqcup B_2$ and $p_1^{\gcd(B_1-B_1)}=1$. Without loss of generality, suppose $\lvert B_1 \rvert \ge 2$, and let $a_i$ and $a_j$ be any two of its elements. So, $p_1^{a_i-a_j}=1$. By Lemma \ref{twotwounity}, $\rk M_A^{2,2} (1,p_1) = 4$, which is a submatrix of $M_A^{2,2,2} (1,p_1,p_2)$. 
\end{proof}

\section{Singularities of sparse discriminants}

The discriminant $D_A\subset\C^A$ is the closure of the discriminant of the projection $p_A$ of $S=\{(f,t)\,|\,f(t)=0\}\subset\C^A\times\C^*$ onto $\C^A$. We study singularities of $D_A$ and $p_A$. 

\begin{theor}\label{mainth}
Assume that $A=\{a_1<a_2<\cdots<a_n\}$ satisfies Convention \ref{conven1}.

1. The components of $\sing D_A$ outside the open torus $\CC^A$ are: 

-- $D_{\min}:=\{c_{a_1}=c_{a_2}=0\}$, provided $a_3-a_2>1$;

-- $D_{\max}:=\{c_{a_n}=c_{a_{n-1}}=0\}$, provided $a_{n-1}-a_{n-2}>1$;

2. A general polynomial $f\in D_{\min}$ has no multiple roots in $\C^*$, the transversal singularity of $D_A$ at this point is an ordinary cusp if $a_1=a_2-1=a_3-3$, and has strictly greater Milnor number otherwise; same for $f\in D_{\max}$.

3. The irreducible set $D_A\cap\CC^A$ splits into the following disjoint parts:

-- The smooth open dense set $U$ of all $f$ having one double root $t_0$ and no other multiple roots; the singularity of $p_A$ at such $(f,t_0)$ is a fold; 

-- A smooth relatively open hypersurface $D_\prec$ of all $f$ having one triple root $t_0$, and no other multiple roots; the singularity of $p_A$ at such $(f,t_0)$ is a cusp;

-- A smooth relatively open hypersurface $D_\times$ of all $f$ having exactly two double roots, and no other multiple roots; the singularity of $D_A$ at such point is a union of two transvesal lines; the multisingularity of $p_A$ at such $f$ is a double fold.

--  the rest $\Sigma$, having codimension greater than 1 in $D_A$.

4. The set $\sing D_A\cap \CC^A$ is contained in $D_\prec\sqcup D_\times\sqcup\Sigma$.
\end{theor}
(One can actually prove the latter is the closure of the former.)

\begin{proof}
    The proof occupies the rest of the section. 
    
    Subsection 3.1 contains the statements about versal deformations that we use in the rest of the proof. In particular, we deduce several results about the discriminant $\nabla$ of a linear system $L$ on a projective variety $X$. 

    Then, we apply these statements to the case $L = \C^A$, $X = \mathbb P^1 \supset \C^*$. Then, $\nabla \cap \CC^A = D_A \cap \CC^A$.
    
    Section 3.2 contains the proof of part 3 of the theorem. We subdivide $D_A \cap \CC^A$ into strata $D_A^{m_1,\dots,m_k}$ according to the multiplicities and numbers of multiple roots. For example, $D_\times = D_A^{1,1} \cap \CC^A$ and $D_\prec = D_A^2 \cap \CC^A$. Proposition \ref{stratadimension} reduces the study of dimensions and transversal types of singularities to computation of the ranks of Vandermonde type matrices in Section 2. In fact, we show that $\C^A$ is a versal deformation of any $f\in D_\times \cup D_\prec$. The necessary ranks are listed in Theorem \ref{mainrank}. The actual proof of the theorem is then contained in in Theorems \ref{dimclosed} and \ref{mainnotinfty}. 

    The stratum $D_A^{1,1,1}$ is treated slightly differently: we also needed to apply basic facts about the dimension of varieties since $\rk M_A^{1,1,1}$ is often less than 3 (see Proposition \ref{oneoneoneabsone}).

    Parts 1 and 2 are proved in Section 3.3. We consider singularities of the curve $W \cap D_A$ where $W\subset \C^A$ is a 2-dimensional affine subspace. The result is derived from the theory of sparse curve singularities \cite{esv}.

    Part 4 trivially follows from the previous parts.\end{proof}

\begin{rem}\label{wildcomponents}
    The variety $D_\prec$ is irreducible. Its explicit rational parameterization was found in \cite{mts}. Interestingly, $D_\prec$ is the projection of another $A$-discriminantal variety that corresponds to polynomials not in one but in two variables.

    We do not know whether the variety $D_\times$ is irreducible under the assumptions we imposed on $A$. However, if $A=(-s,-t,0,t,s)$ for some positive distinct odd integers $s$ and $t$, then $D_\times$ is reducible. One may show that there exist complex numbers $\alpha\ne 1$ such that $\rk M_{A}^{2,2}(1,\alpha) = 3$. Then, for each such $\alpha$, consider the set $D_A^\alpha$ of polynomials having two multiple roots the quotient of which equals $\alpha$. The variety $D_A^\alpha$ is an irreducible component of $D_\times$ of dimension $3$. The remaining polynomials in $D_\times$ form another irreducible component of dimension 3.

    If $\lvert A \rvert = 4$, then $D_\times$ is usually reducible as well since the quotient of two multiple roots of polynomials in $\C^A$ may obtain only a finite number of values.
\end{rem}

\begin{rem}
    One may also show that each irreducible component of the singular locus of $D_A$ has dimension $n-2$. If $f$ belongs to the singular locus and $D_A$ has two branches at $f$, then their intersection has dimension $n-2$ and belongs to the singular locus. If the branch is unique, then $f$ has a root of multiplicity greater than 2, so $f\in \overline{D_\prec}$ and $\dim D_\prec = n-2$.
\end{rem}

\subsection{Preliminaries on versal deformations}

\subsubsection{V-versal deformations and Tjurina algebra}

Consider a germ of a function $f\colon (\C^m,0) \to (\C,0)$.

\begin{defin}
	A germ of analytic function $F(x,\lambda) \colon (\C^m\times \C^n,0) \to (\C,0)$ is called a deformation of $f$ if $F(x,0) = f(x)$ for any $x\in \C^m$. The space $(\C^n,0)$ is called the base of a deformation.
\end{defin}

\begin{defin}
	The function $f\colon \C^m \to \C$ has a singularity at point $p$ if $f(p)=0$ and $df(p)=0$. 

	The gradient ideal $I_f(p)$ of the function $f$ at $p$ is the ideal of the local ring $\mathcal{O}_{p,\C^m}$ generated by the $m$ partial derivatives of $f$. 
	The local algebra $Q_f(p)$ of $f$ is the quotient $\mathcal{O}_{p,\C^m} / I_f(p)$. The dimension $\dim Q_f(p)=\mu_f$ is called the Milnor number of $f$.

    The Tjurina ideal $J_f(p)$ is generated by the partial derivatives and $f$ itself.
	The Tjurina algebra $T_f(p)$ is $\mathcal{O}_{p,\C^m} / J_f(p)$. The dimension $\dim T_f(p)=\tau_f$ is called the Tjurina number of $f$. 

    The Tjurina and Milnor numbers of isolated singularities are finite.
\end{defin}

\begin{defin} \cite[Section 6.5]{avg}
    Germs $f,\tilde f \colon (\C^m,0) \to (\C,0)$ are $V$-equivalent if there exits a germ of a biholomorphic map $g \colon (\C^m,0) \to (\C^m,0)$ and a germ $M\colon (\C,0) \to \C^*$ such that $\tilde f(x) = M(x)f(g(x))$.
\end{defin}

\begin{defin}
    A deformation $\tilde F(x,\lambda) \colon (\C^m \times \C^\alpha,0) \to (\C,0)$ of the function $f$ is called $V$-versal if any deformation $F(x,\lambda) \colon (\C^m \times \C^\beta,0) \to (\C,0)$ can be induced from $\tilde F$, that is there exist germs of maps $M \colon (\C^m \times \C^\beta,0) \to \C^*$, $g \colon (\C^m \times \C^\beta,0) \to (\C^m,0)$, and $\varphi \colon (\C^\beta,0) \to (\C^\alpha,0)$ such that $F(x,\lambda) = M(x,\lambda) \tilde F(g(x,\lambda),\varphi(\lambda))$ and $g(x,0) \equiv x$, and the number $M(0,0)$ is nonzero. 

    The deformation $\tilde F$ is called miniversal if $\alpha$ equals the Tjurina number $\tau$ of $f$. Such a deformation always exists and is essentially unique.
\end{defin}

The presence of the map $M$ allows us to regard $f$ not only as a function, but also as a section of a linear bundle: indeed, the quotient of two sections is a locally non-vanishing holomorphic function.

Consider a germ of an isolated singularity $f\colon (\C^m,0) \to (\C,0)$ and its deformation $F\colon (\C^m\times\C^k,0) \to (\C,0)$, $F(x,0)=f(x)$. By $\dot F_i$ denote the germs of functions ${\partial F(x,\lambda) \over \partial \lambda_i} \vert_{\lambda=0} \in \mathcal O_{p,\C^m}$. We call them velocity vectors.

By $r$ denote the dimension of the image of the natural projection $\langle \dot F_1, \dots, \dot F_k \rangle_\C \to T_f$ from the linear span of the velocity vectors in $\mathcal O_{p,\C^m}$ to the Tjurina algebra.

It is known \cite[Section 8.3]{avg} that $r = \dim T_f$ if and only if the deformation $F$ is $V$-versal. In this case, $F$ is called an infinitesimally versal deformation of $f$

\begin{lemma}\label{deformrankV} \cite[Section 8.5]{avg}
	
	Let $F\colon (\C^m\times\C^n,0) \to (\C,0)$ be any deformation of $f$, and $\tilde F(x)\colon (\C^m\times\C^\tau,0) \to (\C,0)$ be a $V$-miniversal deformation of $f$. As usual, the deformation $F$ is induced from $\tilde F$: \begin{equation} \label{induceV}
		F(x,\lambda) = M(x,\lambda) \tilde F(g(x,\lambda),\varphi(\lambda)).
	\end{equation}
	
	Then the rank of the map germ $\varphi \colon (\C^n,0) \to (\C^\tau,0)$ at 0 equals $r$.
	
\end{lemma}
\begin{proof}
	Let us compute the differential

	\[{dF \over d\lambda}  = {dM \over d \lambda} \tilde F + M {\partial \tilde F \over \partial g}{dg \over d\lambda}  + M{\partial \tilde F \over \partial \varphi}{d\varphi \over d\lambda}.\]
	
	Using the identity $\tilde F(x,0) = f(x)$ we obtain that 
	
	\[{dF \over d\lambda} \vert_{\lambda=0} = f {dM \over d \lambda}\vert_{\lambda=0}  + M{\partial f \over \partial x}{dg \over d\lambda} \vert_{\lambda=0} + M{\partial \tilde F \over \partial \varphi}{d\varphi \over d\lambda} \vert_{\lambda=0}.\]
	
	Here, ${dF \over d\lambda}\vert_{\lambda=0}$ is the row of velocity vectors $(\dot F_1, \dots, \dot F_n) \in (\mathcal{O}_{0,\C^m})^n$, ${dM \over d \lambda}$ is a row consisting of elements in $\mathcal{O}_{\C^m}$, ${\partial f \over \partial x}$ is the row of partial derivatives $(f'_{x_1},\dots,f'_{x_m}) \in (J_f)^m$, ${dg \over d\lambda} \vert_{\lambda = 0}$ is a matrix in $Mat_{m \times n}(\mathcal{O}_{0,\C^m})$, ${\partial \tilde F \over \partial \varphi}$ is the row of velocity vectors $(\dot {\tilde F}_1,\dots, \dot {\tilde F}_\tau)$ and ${d \varphi \over d \lambda}\vert_{\lambda=0}\in Mat_{\tau \times n}(\C)$ is the differential of the map $\varphi$.
	
	Note that the row ${\partial f \over \partial x}{dg \over d\lambda} \vert_{\lambda=0}$ belongs to $(J_f)^n$ since its elements are $\mathcal{O}_{0,\C^m}$-linear combinations of partial derivatives of $f$. Moreover, $f {dM \over d \lambda}\vert_{\lambda=0} \in (J_f)^n$.
	
	By $[z]$ denote the class of the function $z \in \mathcal{O}_{\C^m}$ in $T_f$. It now follows that
	
	\[([\dot F_1], \dots, [\dot F_n]) = ([\dot {\tilde F}_1],\dots, [\dot {\tilde F}_\tau]) {d\varphi \over d\lambda} \vert_{\lambda=0}.\]
	
	Since $\tilde F$ is a miniversal deformation, the vectors $[\dot {\tilde F}_1],\dots, [\dot {\tilde F}_\tau]$ form a basis of the vector space $T_f$. The columns of the matrix ${d\varphi \over d\lambda} \vert_{\lambda=0}$ encode the representation of the velocity vectors of $F$ as linear combinations of these basis vectors. Hence, the rank of this matrix is equal to the dimension $r$ of the linear span of the vectors $[\dot F_i]$.
\end{proof}

\begin{lemma} \label{zeropreimage}
	Assume that $F$ is a miniversal deformation of $f\colon (\C^n,0) \to (\C,0)$. Then there exists a neighborhood $U\subset \C^n$ with the following property. For any $\lambda_0\in U$ the function $F(x,\lambda_0)$ is $V$-equivalent to $f$ if and only if $\lambda_0=0$.
\end{lemma}

\begin{proof} 
We will prove this lemma in the case in which we are particularly interested: $n=1$. Since a function of one variable $f$ is quasihomogeneous, the Tjurina and Milnor algebras coincide and thus the $V$-miniversal deformation coincides with $R$-miniversal. It is known \cite{gabrielov} that $F(x,\lambda_0)$ is $R$-equivalent to $f$ if and only if $\lambda_0=0$. 

In general, the proof is obtained by reduction the problem to finite-dimensional spaces of functions by Tougeron's theorem, and regarding the action of the group of local diffeomorphisms as an action of a finite-dimensional algebraic group.
\end{proof}

\subsubsection{Bifurcation diagrams and transversal types of discriminants}

Consider a germ of an isolated singularity $f \colon (\C^n,0) \to (\C,0)$ and its miniversal deformation $F(z,\lambda) \colon (\C^n \times \C^\tau,0) \to (\C,0)$.

\begin{defin}
	The $V$-bifurcation diagram of the singularity of $f$ is the set of functions in the base space $\C^\tau$ of the $V$-miniversal deformation of $f$ that have critical value 0.
\end{defin}

\begin{exa} \label{swallowtail}
	A function $f\colon (\C,0) \to (\C,0)$ has singularity of type $A_n$ if it has a root of multiplicity $n+1$ at zero. The bifurcation diagram of $A_n$ is called a generalized swallowtail. We denote its germ at 0 by $\Sigma_n$. In particular, $\Sigma_1 \cong \{0\} \subset \C$,   $\Sigma_2 \cong \{x^2-y^3=0\} \subset \C^2$. 
\end{exa}

Let $V\subset H^0(L)$ be a linear system in a proper smooth variety $X$. We will often regard elements of $V$ as rational functions that do not vanish in a neighborhood of several points: for any $s_0 \in H^0(L)$ there exists an embedding $H^0(L) \to \C(X)$, $s\mapsto s/s_0$.

\begin{defin}
	The discriminant $\nabla_V$ is the set of all sections in $V$ that define a divisor that has a singular point.
\end{defin}

\begin{defin}
	The tautological variety of the linear system $V$ is the variety $\nabla'_V \subset V \times X$ that consists of the pairs $(f,x)$ such that $f$ defines a divisor $D$ such that $x$ is a singular point of $D$.
\end{defin}

\begin{lemma} \label{neighborhoods}
	Consider $f\in \nabla$ corresponding to a divisor with singular points $p_1, \dots, p_k$ and no other singular points. For each $i$, consider an metric neighborhood $U_i$ of $p_i$. Then there exists a metric neighborhood $V_{0}\subset V$ of $f$ such that each divisor corresponding to an element of $\nabla \cap V_{0}$ has a singular point in one of the sets $U_i$.

    As a corollary, for each branch of $\nabla$ at $f$ there exists $i$ such that the functions on that branch have a singular point in $U_i$.
\end{lemma}

\begin{proof}
Assume that such a neighborhood $V_{0}$ does not exist. Then there exists a sequence $f_j$ of elements of $\nabla$ converging to $f$ such that the divisor of $f_j$ has no singularities in the union of the sets $U_i$. For each $j$, choose a singular point $x_j\in X$ of the divisor of $f_j$. Since $X$ is compact, we may assume that the sequence $(x_j)$ converges to some point $x$. The sequence $(f_j,x_j)$ of elements of $\nabla'$ converges to $(f,x)$. But $\nabla'$ is closed in $X \times V$ since it is defined by system of algebraic equations $f = df = 0$. Thus, $x$ is a singular point of $f$. Hence, $x \in U_i$ which is impossible. 
\end{proof}

\begin{theor} \label{universe}
	Choose an element $f$ of $V$ such that has a singularity at exactly one point of $X$.
	
	By $W$ denote the set of all functions in $V$ that have exactly one singularity at some point that is $V$-equivalent to the singularity of $f$ at $p$. 
	Then $W$ if a locally closed analytic subset of $V$. 
	
	Moreover, by $r$ denote the dimension of the image of $V$ in the local algebra $T_f$ under the quotient map. 
	Then the local dimension of $W$ at point $f$ satisfies the inequality $\dim_f W \le \dim V - r$.
\end{theor}

\begin{proof}
	 The function $f$ has exactly one singular point at $p$. By $\pi_f$ denote the natural map of $V$ to $T_f(p)$. By $\varphi_f \colon (V,f) \to \C^\tau$ denote the natural map to the versal deformation of the singularity $f$. Then by lemma \ref{zeropreimage}, there is a neighborhood $U$ of $f$ in $V$ such that $U\cap W = \varphi_f^{-1}(0)$. Certainly, $\varphi_f^{-1}(0)$ is an analytic set closed in $U$.
	 
	  Note that in the previous paragraph we could replace $f$ by  any function that belongs to $W$ (and replace $p$ by the singular point of this function). Cover $W$ by the open analytic subsets $U$ corresponding to such functions. Hence, the first claim of the theorem holds.
	  
	  Now, since the rank of the map $\varphi_f$ is $r$ by lemma \ref{deformrankV}, $\dim(U\cap W) = \dim(\varphi_f^{-1}(0)) \le \dim V-r$. 
	\end{proof}

\begin{theor} \label{multiverse}
	 Choose an element $f$ of $V$ such that has singularities at points $p_1,\dots,p_k$ of $X$ and no other singularities.
	
	By $W$ denote the set of all functions in $V$ that have singularities at some points $q_1, \dots, q_k$ equivalent to the singularities of $f$ at $p_1, \dots, p_k$ respectively and no other singularities.
	Then $W$ if a locally closed analytic subset of $V$. 
	
	Moreover, by $r$ denote the dimension of the image of $V$ in the direct sum $T_f(p_1) \oplus \dots \oplus T_f(p_k)$ under the direct sum of the quotient maps. 
	Then $\dim_f W \le \dim V - r$.
\end{theor}

\begin{proof}
	Similar to the proof of Theorem \ref{universe}.
\end{proof}

Let $W$ be a smooth $m'$-dimensional manifold contained in $\nabla$.

\begin{defin}
	Let $p$ be a point of $W$ and $Y$ be a germ of analytic subset of $(\C^{m-m'},0)$. The pair $(\nabla,W)$ has transversal type $Y$ at $p$ if the germ at $p$ any smooth hypersurface slice of $\nabla$ transversal to $W$ is isomorphic to $Y$. That is, for any manifold $H$ with $\dim H = m-m'$ passing through $p$ transversally to $W$, the germ $(\nabla \cap H,p) \subset H$ is the image of $(Y,0) \subset \C^{m-m'}$ under some biholomorphic map $(Y,0) \to (H,p)$.
\end{defin}

\begin{theor} \cite[Corollary 21.4.4]{avg} \label{transversalfullrank}
	Under conditions of Theorem \ref{multiverse}, assume moreover that $r=\sum_i\tau_f(p_i)$. By $\C^{\tau_i}$ denote the base space of the miniversal deformation of the singularity of $f$ at $p_i$. By $B_i\subset \C^{\tau_i}$ denote the bifurcation diagram of $f$ at $p_i$.
	
	Then $\dim_f W = \dim V - r$, the point $f$ of the variety $W$ is smooth and the transversal type of $\nabla$ at $p$ with respect to $W$ is 
	\[B=\{(\lambda_1,\dots, \lambda_k)\in \C^{\tau_1} \times \C^{\tau_2} \times \dots \times \C^{\tau_k} \mid \lambda_i \in B_i \text{ for some } i\}\]
\end{theor}

\begin{proof}
	Consider the canonical map to the versal deformation of $f$: $\varphi \colon (V,f) \to (\C^{\tau_1} \times \dots \times \C^{\tau_k},0)$. By Lemma \ref{deformrankV}, $\varphi$ is of full rank $\sum \tau_i$. 
	
	Consider a small neighborhood $U$ of $f$. Each element of $\nabla \cap U$ corresponds to a divisor on $X$ that has a singularity in a small neighborhood of one of the points $p_i$. Thus, by lemma \ref{neighborhoods}, the $\varphi$-image of a point of $\nabla\cap U$ is contained in $B_i$ for some $i$. Hence, $\nabla \cap U = \varphi^{-1}(B)$.
	
	By lemma \ref{zeropreimage}, $W\cap U = \varphi^{-1}(0)$. Since $\varphi$ is a local diffeomorphism, $W \cap U$ is smooth. For any smooth subvariety $V'$ of $V$ passing through $f$ transversally to $W$ we find that $\varphi|_{V'} \colon V' \to \C^{\tau_1} \times \dots \times \C^{\tau_k}$ is a local diffeomorphism, and thus $\nabla \cap V'$ is diffeomorphic to $B$.
\end{proof}

\begin{proof}
    The image of $V$ under $\varphi$ is a germ of analytic variety of dimension at least 2. Thus, $\varphi(V\setminus \{f\})$ intersects each of the bifurcation diagrams $B_i\times\C^{\sum_j(\tau_j) - \tau_i}$ which are hypersurfaces. Note that $\varphi(\nabla)$ is not contained in any of the bifurcation diagrams (otherwise, a generic element of $V$ has a singular point). The $\varphi$-preimage of each of the bifurcation diagrams is a branch of $\nabla$ at $f$. Thus, $f$ is a singular point if $k \ge 2$
    
    Now, assume $k=1$. Then $\tau_1\ge2$. The description of bifurcation diagrams of simple singularities is well-known, and the equality $\dim B_1^{sing} = \dim B_1 - 1$ holds for them if $\tau_1\ge 2$. 
\end{proof}

\subsubsection{Folds and cusps}
\begin{utver}(see \cite[Section 9.3]{avg}) \label{whitneycusp}
    Let $F(x,\lambda) \colon (\C_x \times \C_{\lambda_0,\dots,\lambda_{m-1}}^m,(x_0,0)) \to (\C,0)$ be a versal deformation at $x_0$ of an analytic function $f(x)$ having root of multiplicity $k\ge 2$ at the point $x_0$. Denote $\nabla'=\{(x,\lambda)\in \C\times \C^m \mid F(x,\lambda)=0 \}$. Then the projection $p \colon \nabla' \to \C^m$ has a Thom-Boardman singularity of type  $S_{k-1}$ singularity at the point $(x_0,f)$.

    In particular, if $k=2$, then $p$ has a fold and if $k=3$, then $p$ has a cusp at  $(x_0,f)$. 
\end{utver}

\begin{proof}
    We may choose  coordinates on $\C_x$ so that $x_0=0$ and $f(x)=x^k$.
    
    We will write $(\lambda_0,\dots,\lambda_{m-1})=\lambda$ and $(\lambda_1,\dots,\lambda_{m-1})=\vec \lambda$.
    By $\tilde F(x,\mu_0,\mu_1,\dots,\mu_{k-2})=\mu_0 + \mu_1 x + \dots + \mu_{k-2}x^{k-2} + x^k$ denote the miniversal deformation of $f$.

    Then $F(x,\lambda) = M(x,\lambda) \tilde F(g(x,\lambda),\mu_0(\lambda),\dots, \mu_{k-2}(\lambda))$. Since $F$ is a versal deformation, the map $\lambda \mapsto (\mu_0(\lambda),\dots,\mu_{k-2}(\lambda))$ is of rank $k-1$ (Lemma \ref{deformrankV}). Thus, we may choose the coordinates on $\C^m_\lambda$ so that $\mu_i(\lambda)=\lambda_i$ for all $i$.

    Now, the variety $\nabla'$ is defined by equation $\lambda_0 + \lambda_1 g(x,\lambda) + \dots + \lambda_{k-2} g(x,\lambda)^{k-2}+ g(x,\lambda)^k=0$.

    Since $g(x,0)=x$, $(y=g(x,\lambda),\lambda_0,\dots,\lambda_{m-1})$ is a coordinate system on $\C\times\C^m$; denote $x=g^{-1}(y,\lambda).$ Thus, $\nabla'=\{(y,\lambda)\mid \lambda_0 + \dots + \lambda_{k-2} y^{k-2} + y^k=0\}$.

    The variables $y,\lambda_1,\dots,\lambda_{m-1}$ define a coordinate chart on $\nabla'$. The map $p$ on this chart is $(y,\vec \lambda) \mapsto (-\lambda_1 y - \lambda_2 y^2 -\dots -\lambda_{k-2} y^{k-2} -y^k, \vec \lambda)$. This is one of the standard coordinate forms of the Thom-Boardman singularity.
\end{proof}

\subsection{Dimensions of discriminant strata}

By $\C^{A-}$ denote the set of all polynomials $\sum_{a\in A} c_a t^a$, $c_{a_1}\ne 0$, $c_{a_n}\ne 0$.
\begin{utver}
    Consider a polynomial $f = \sum_{a\in A} c_a t^a$, $c_{a_1}\ne 0$, $c_{a_n}\ne 0$. Then $f \in D_A$ if and only if $f$ has a nonzero root of multiplicity greater than one.
\end{utver}

\begin{proof}
    The "if" part follows from the definition of $D_A$. 

    Now, consider $f\in \C^{A-} \cap D_A$. Assume that $f$ all the nonzero roots of $f$ are simple. Then there exists a sequence $f_n$ of elements of $D_A^\circ$ that converges to $f$. 
    
    For each $n$, choose a multiple root $x_n \in \C \setminus 0$ of $f_n$. The sequence $x_n$ has a subsequence that converges to a point $x_0 \in \mathbb P^1$. If $x_0 \in \C \setminus 0$, then $x_0$ is a multiple root of $f$. 
\end{proof}

\begin{defin}
    Assume that $f\in \C^{A-} \cap D_A$. Let $\{p_1, \dots, p_k \}\in \C^*$ be the set of roots of $f$ of multiplicity greater than one. By $m_i+1$ denote the multiplicity of the root $p_i$. By $T_f(p_i) = \C[t^{\pm1}]_{p_i}/(f(t),f'(t)) \cong \C^{m_i}$ denote the Tjurina algebra of the singularity of $f$ at $p_i$. Here, $\C[t^{\pm1}]_{p_i}$ is the localization of the Laurent polynomial ring at the closed point $p_i$. By $T_f = \C[t^{\pm1}]/(f,f')\cong \bigoplus_i T_f(p_i)$ denote the global Tjurina algebra of $f$. 
\end{defin}

\begin{defin}
    Let $p_1,\dots p_k$ be a collection of distinct points in $\C \setminus 0$. 
    
    By $D_A^{m_1,\dots,m_k +}(p_1,\dots,p_k)$ denote the vector subspace of $\C^A$ consisting of the functions having roots of multiplicity at least $m_i+1$ at $p_i$ for each $i$. By $D_A^{m_1,\dots,m_k}(p_1,\dots,p_k)$ denote the subset of $\C^{A-}$ consisting of the functions having roots of multiplicity exactly $m_i+1$ at each $p_i$ and no other roots of multiplicity greater than 1.

\end{defin}

\begin{defin}
    By $\CC^k\setminus diag \subset \CC^k$ denote the set of configurations of distinct points. The multisingularity stratum is the following variety.

    \[D_A^{m_1,\dots,m_k} = \bigcup_{(p_1,\dots,p_k) \in \CC^k \setminus diag} D_A^{m_1,\dots,m_k}(p_1,\dots,p_k). \]

It will also by convenient to consider the following two sets.
\[D_A^{m_1,\dots,m_k +} = \bigcup_{(p_1,\dots,p_k) \in \CC^k \setminus diag} D_A^{m_1,\dots,m_k +}(p_1,\dots,p_k)\]

  If $U\subset \CC^k \setminus diag$ is any subset, then \[D_A^{m_1,\dots,m_k}(U) = \bigcup_{(p_1,\dots,p_k) \in U} D_A^{m_1,\dots,m_k}(p_1,\dots,p_k) \]
\end{defin}

  \begin{observ} \label{fiberdim}
    The set $D_A^{m_1,\dots,m_k+}(p_1,\dots,p_k)$ is a vector subspace of $\C^A$. Its dimension is $n-\rk M_A^{m_1+1,m_2+1,\dots,m_k+1}(p_1,\dots,p_k)$. 
\end{observ}

  \begin{defin}
  We define a partial order on the set of all finite sequences of integers as follows. For sequences  $S = (S_i)_{i=1\dots k}$ and $\tilde S = (\tilde S_i)_{i=1\dots \tilde k}$ consider their rearrangements $S^\downarrow$ and $\tilde S^\downarrow$ in nonincreasing order.  Then $S \succeq \tilde S$ if $k\ge \tilde k$ and $S^\downarrow_i\ge \tilde S^\downarrow_i$ for $i\le \tilde k$.
  \end{defin}

  \begin{observ} \label{rankineq}
      Assume that $\tilde S \succeq S$ are two nonincreasing sequences of lenghthes $\tilde k$ and $k$ respectively. Let $p_1,\dots,p_{\tilde k}\in \C^*$ be distinct points. Then $\rk M_A^{\tilde S}(p_1,\dots,p_{\tilde k}) \ge \rk M_A^S(p_1,\dots,p_k)$.
  \end{observ}

  \begin{observ} \label{stratumclosure}
    Assume that $m_1\ge m_2 \ge \dots \ge m_k$. Then $D_A^{m_1,\dots,m_k+} \cap \C^{A-}$ is the union of all sets $D_A^{\tilde m_1, \tilde m_2\dots, \tilde m_{\tilde k}}$ for which $(\tilde m_1,\dots,\tilde m_{\tilde k}) \succeq (m_1,\dots,m_{k})$.

    For any $A$, this union consists of finite number of sets.
  \end{observ}

  \begin{utver} \label{stratadimension}
      The set $D_A^{m_1,\dots,m_k}$ is analytic and locally closed. Assume that $f\in D_A^{m_1,\dots,m_k} (p_1,\dots,p_k)$ and denote $\rk M_A^{m_1,\dots,m_k}(p_1,\dots,p_k) = r$. Then 
      
      A) The inequality $\dim_f D_A^{m_1,\dots,m_k} \le n-r$ holds.

      B) Assume that $r = m_1 + \dots + m_k$. Then $f$ is a smooth point of $D_A^{m_1,\dots,m_k}$ and $\dim_f D_A^{m_1,\dots,m_k} = n-r$. 
      
      By $\Sigma'_{m_i}$ denote the preimage of the bifurcation diagram $\Sigma_{m_i} \subset \C^{m_i}$ (see Example \ref{swallowtail}) under the projection $\C^{m_1} \times \C^{m_2} \dots \times \C^{m_k} \to \C^{m_i}$. The transversal type of the pair $(D_A,D_A^{m_1,\dots,m_k})$ at $f$ is \[\Sigma'_{m_1} \cup \dots \cup \Sigma'_{m_k} \subset \C^{m_1+\dots + m_k}.\]

      C) Let $U$ be an open subset of $\C^*$. Assume that for any $p_1,\dots, p_k \in \C^*$ the inequality $\rk M_A^{m_1,\dots,m_k}(p_1,\dots,p_k) \ge r$ holds. Then $\dim D_A^{m_1,\dots,m_k}(U) \le n-r$ and $\dim D_A^{m_1,\dots,m_k+}(U) \le n-r$.
  \end{utver}

  \begin{proof}
      The torus $\C^*$ is naturally embedded in $\Pr^1$. The vector space $\C^A$ can be considered as a linear system on $\Pr^1$: this is a subspace of $\Gamma(\mathcal O (a_n-a_1),\Pr^1)$.

      By Observation \ref{vandermonderank}, $r$ is the rank of the natural map to the Tjurina algebra of $f$:

      \begin{equation*}
          \pi_f \colon \C^A \to \C^{m_1} \oplus \dots \oplus \C^{m_k}, \; \rk \pi_f = r
      \end{equation*}.

      Now, we apply Theorem \ref{multiverse} to $f$. The set $W$ of this theorem is a union $D_A^{m_1,\dots,m_k} \cup D_A^\infty$ where $D_A^\infty$ is a certain subset of $\C^A \setminus \C^{A-}$. We obtain $D_A^{m_1,\dots,m_k}=W\setminus (\C^A \setminus \C^{A-})$, so $D_A^{m_1,\dots,m_k}$ is a locally closed analytic set. Since $f \notin D_A^\infty$, $\dim_f D_A^{m_1,\dots,m_k} = \dim_f W \le n-r$.

  Now, part B of the proposition follows from Theorem \ref{transversalfullrank}.

  The inequality $\dim D_A^{m_1,\dots,m_k}(U) \le n-r$ of part C follows from the fact that for any $f\in D_A^{m_1,\dots,m_k}(U)$ the local dimension $\dim_f D_A^{m_1,\dots,m_k}(U) \le n-r$.

  The inequality $\dim D_A^{m_1,\dots,m_k+}(U) \le n-r$ now follows from Observation \ref{stratumclosure}.
  \end{proof}

\begin{theor} \label{dimclosed} Assume that $n\ge 6$ and that $A$ is doubly reduced. Then:
\[\dim D_A^{2+} \le n-2; \; \; \; \dim D_A^{1,1+} \le n-2;\]
    \[\dim D_A^{3+} \le n-3; \; \; \; \dim D_A^{2,1+} \le n-3; \; \; \; \dim D_A^{1,1,1+} \le n-3.\]
  \end{theor}
\begin{proof}
    The first four inequalities are already essentially proved. Substitute the rank estimations of Theorem \ref{mainrank} (1), (11), (1), (21) in Proposition \ref{stratadimension} C.

    For example, $D_A^{2,1+} = \bigcup_{\mathbf m \succeq (2,1)} D_A^{\mathbf m}$. By Theorem \ref{mainrank}  (21), $\rk M_A^{2,1} \ge 3$. Thus, by observation \ref{rankineq}, for any $\mathbf m \succeq (2,1)$, $(p,q) \in \CC^2$ the inequality $\rk M_A^{\mathbf m}(p,q) \ge 3$ holds. Hence, Proposition \ref{stratadimension} C (with $U = \CC^2 \setminus diag$), implies the inequality $ \dim D_A^{\mathbf m} \le n-3$.

    Now, let us prove the fifth inequality. Consider the set $X_A$ of Theorem \ref{mainrank}. By $X_A'$ denote the variety $\{(p_1,p_2,p_3) \in \CC^3\mid (p_1^{-1}p_2,p_1^{-1}p_3)\in X_A, p_1\ne p_2 \ne p_3 \ne p_1\}$. Let $U = \CC^3 \setminus X_A'$. By definition of $X_A$ and Observation \ref{matrixtorusaction}, for any $(p_1,p_2,p_3) \in U$, $\rk M_A^{1,1,1} (p_1,p_2,p_3) =3$. Thus, by Proposition \ref{stratadimension}, $\dim D_A^{1,1,1+}(U) \le n-3$.

    Note that $D_A^{1,1,1+} = D_A^{1,1,1+}(U) \cup D_A^{1,1,1+}(X_A')$. Thus, it is enough to show that $\dim D_A^{1,1,1+}(X_A') \le n-3.$

    By Theorem \ref{mainrank} (111), $X_A$ is a finite set, so $X_A'$ is a curve. Consider the variety \[Y = \{(f,p_1,p_2,p_3) \in \C^A  \times X_A' \mid f \text{ has roots of multiplicity at least 2 at $p_1,p_2$ and $p_3$}\}\] By $\pi_1$ and $\pi_2$ denote the projections of $Y$ onto $\C^A$ and $X_A'$ respectively. Clearly, $\pi_1(Y)\supset D_A^{1,1,1+}(X_A')$, so $\dim D_A^{1,1,1+}(X_A') \le \dim Y$.

    Now, for any $(p_1,p_2,p_3)\in X_A'$, 
    
    \begin{align*}
        &\dim \pi_2^{-1}(p_1,p_2,p_3)  = \dim D_A^{1,1,1+}(p_1,p_2,p_3) =  [\text{Observ. \ref{fiberdim}}] = \\ &= n-\rk M_A^{2,2,2}(p_1,p_2,p_3) = [\text{Observ. \ref{matrixtorusaction}}] =  n-\rk M_A^{2,2,2}(1,p_1^{-1}p_2,p_1^{-1}p_3) \le \\&\le [\text{Thm \ref{mainrank} (222)}] \le n-4.
    \end{align*}
    Hence, the dimension of a generic fiber of $\pi_2 \colon Y \to X_A'$ is at most $n-4$. Hence, $\dim Y \le \dim X_A' + n-4 = n-3$.
\end{proof}

\begin{utver} \label{dimclosedfrombelow} For any $m_1,\dots,m_k$ the inequality
	$\dim D_A^{m_1,\dots,m_k+} \ge n-m_1-\dots-m_k$ holds.
\end{utver}

\begin{proof}
	Consider the variety \[Y = \{(f,p_1,\dots,p_k) \in D_A^{m_1,\dots,m_k+} \times \CC^k \mid f\in D_A^{m_1,\dots,m_k+}(p_1,\dots,p_k) \}.\] The projection $\pi_1 \colon Y \to D_A^{m_1,\dots,m_k+}$ is a finite map, so $\dim D_A^{m_1,\dots,m_k+}=\dim Y$. The dimension of each fiber of $\pi_2 \colon Y \to \CC^k$ is at least $n-(m_1+\dots+m_k) -k$ by Lemma \ref{fiberdim}, so $\dim Y \ge n-(m_1+\dots+m_k) -k + \dim \CC^k = n-(m_1+\dots+m_k)$.
\end{proof}

\begin{utver} \label{mainnotinfty}
	(1) The points of $D_A^1$ are smooth in $D_A$. $\dim D_A^1=n-1$.
	
	(2) The set $D_A^2$ is a smooth variety and $\dim D_A^2=n-2$. The pair $(D_A,D_A^2)$ has transversal type $x^2-y^3=0$ at any point.
	
	(3) The set $D_A^{1,1}$ is a smooth variety and  $\dim D_A^{1,1} = n-2$. The pair $(D_A,D_A^{1,1})$ has transversal type $x^2-y^2=0$ at any point.
	
	(4) $\dim ((D_A\cap \C^{A-})\setminus D_A^1 \setminus D_A^2 \setminus D_A^{1,1}) \le n-3$.

    (5) The singularity of $p \colon S \to \C^A$ at any point of $D_A^1$, $D_A^2$, and $D_A^{1,1}$ is, respectively, a fold, cusp, and double fold.
\end{utver}

\begin{proof}
	By Observation \ref{stratumclosure}, $D_A^1 = D_A^{1+} \setminus D_A^{2+} \setminus D_A^{1,1+}$, $D_A^2 = D_A^{2+} \setminus D_A^{3+} \setminus D_A^{2,1+}$, $D_A^{1,1} = D_A^{1,1+} \setminus D_A^{2,1+} \setminus D_A^{1,1,1+}$, $(D_A\cap \C^{A-}) \setminus D_A^2 \setminus D_A^{1,1} \subset D_A^{3+} \cup D_A^{2,1+} \cup D_A^{1,1,1+}$. Applying the dimension estimations of Theorem \ref{dimclosed} and Proposition \ref{dimclosedfrombelow} we compute the desired dimensions. 
	
	The statements about smoothness and transversal types are special cases of Proposition \ref{stratadimension} (using the rank estimations of Theorem \ref{mainrank} (1) for $m=2$ and (11)).

    Part 5 of the theorem for folds and cusps follows directly from Proposition \ref{whitneycusp}. If $f\in D_A^{1,1}$, then $p^{-1}(f)$ has two singular points. The singularities at the points are folds, and the neighborhood of $f$ in $D_A$ is homeomorphic to the union of two transversal hyperplanes. Hence, the singularity of $p$ at $f$ is a double fold. 
\end{proof}

\subsection{Singularities at infinity}

\begin{lemma} \label{rationalsections}
    1. Consider a 2-dimensional affine subspace $W\subset \C^A$ defined by 
    \begin{equation} \label{subspaceparam0}
		W = \{(\alpha_1 \lambda_1 + \beta_1 \lambda_2 + \gamma_1,\dots,\alpha_n \lambda_1 + \beta_n \lambda_2 + \gamma_n)\in \C^A \mid \lambda_1,\lambda_2 \in \C\}.
	\end{equation}

    The variables $\lambda_1,\lambda_2$ define a coordinate system on $W$. 

    The intersection $W \cap D_A$ contains the rational curve with the following parameterization (where $i$ ranges from $1$ to $n$ in the sums).

    \begin{equation} \label{sectionparam1}
		\begin{split}
			\lambda_1(x) &= \frac{ \sum \beta_i x^{a_i} \sum a_i \gamma_i x^{a_i}-  \sum a_i \beta_i x^{a_i} \sum \gamma_i x^{a_i}}{ \sum \alpha_i x^{a_i} \sum a_i \beta_i x^{a_i}  - \sum a_i \alpha_i x^{a_i}   \sum \beta_i x^{a_i} };\\
			\lambda_2(x) &=\frac{ -\sum \alpha_i x^{a_i} \sum a_i \gamma_i x^{a_i}+  \sum a_i \alpha_i x^{a_i} \sum \gamma_i x^{a_i}}{ \sum \alpha_i x^{a_i} \sum a_i \beta_i x^{a_i}  - \sum a_i \alpha_i x^{a_i}   \sum \beta_i x^{a_i} }.
		\end{split}
	\end{equation}

    2. Assume $\gamma_1=\gamma_2=0, \alpha_1=1, \beta_1=0, \alpha_2=0, \beta_2=1, \gamma_i \ne 0 \text{ if } i\ge 3$. Then
    \begin{equation} \label{sectionparam2}
		\begin{split}
			\lambda_1(x) &=   \frac{(a_3-a_2)\gamma_3}{a_2-a_1}x^{a_3-a_1} + o(x^{a_3-a_1}); \\
			\lambda_2(x) &=  - \frac{(a_3-a_1)\gamma_3}{a_2-a_1}x^{a_3-a_2} + o(x^{a_3-a_2}). 
		\end{split}
	\end{equation}

    3. Under the same assumption, let $\kappa>1$ be a common divisor of $a_3-a_1$ and $a_3-a_2$. By $N$ denote the smallest integer such that $a_N-a_1$ is not divisible by $\kappa$.

    Then, the monomials with the smallest exponents not divisible by $\kappa$ in the Taylor expansions of $\lambda_1(x)$ and $\lambda_2(x)$ are $\frac{\gamma_N(a_N-a_2)}{a_2-a_1}x^{a_N-a_1}$ and $\frac{\gamma_N(a_1-a_N)}{a_2-a_1}x^{a_N-a_2}$ respectively.
\end{lemma}
\begin{proof}
    1. By definition of $D_A$, the intersection $D_A \cap W$ is the closure set of all points $(\lambda_1,\lambda_2) \in W$ such that the following two conditions hold for some $x \in \C \setminus 0$:
	
	\begin{equation*} 
		\begin{cases}
		 \sum_{i= 1}^n (\alpha_i \lambda_1 + \beta_i \lambda_2 + \gamma_i) x^{a_i}=0;\\
		  \sum_{i= 1}^n a_i(\alpha_i \lambda_1 + \beta_i \lambda_2 + \gamma_i) x^{a_i}=0.
		\end{cases}
	\end{equation*}

    Solving this system of linear equations  for $\lambda_1,\lambda_2$ we obtain the parameterization \eqref{sectionparam1}.

    2. The lowest degree monomials in the numerator and denominator of $\lambda_1$ are $(a_3-a_2) \gamma_3 x^{a_2+a_3}$ and $(a_2-a_1) x^{a_1+a_2}$ respectively. Their ratio is the lowest degree coefficient of the expansion of $\lambda_1$. The first term of $\lambda_2$ is computed in a similar way.

    3. The numerator of $\lambda_1$ is $x^{a_2+a_3}\left( \sum_{i,j=1}^n \beta_i \gamma_j (a_j-a_i) x^{a_i-a_2} x^{a_j-a_3} \right) = x^{a_2 + a_3} P(x)$. The terms of $P(x)$ have nonnegative degrees since $\beta_1=0,\gamma_1=\gamma_2=0$. It is easy to check that the first term of $P(x)$ of a degree not divisible by $\kappa$ is $\gamma_N (a_N-a_2) x^{a_N-a_3}$.

    Similarly, the denominator of $\lambda_1$ is $x^{a_1+a_2} \left((a_2-a_1) + Q(x) \right)$ where $Q(x)$ is a polynomial such that $Q(0)=0$. The first term of $Q(x)$ of a degree not divisible by $\kappa$ is $\beta_N (a_N-a_1) x^{a_N-a_2}$. 

    Hence, $\lambda_1(x) = x^{a_3-a_1} \frac{P(x)}{(a_2-a_1) + Q(x)} = \left(\frac{x^{a_3-a_1}P(x)}{a_2-a_1}\right)\left(1 -  \frac{Q(x)}{a_2-a_1} + \frac{Q(x)^2}{(a_2-a_1)^2} - \dots \right)$. Since the lowest degrees of terms $P$ and $Q$ not divisible by $\kappa$ are  $a_N-a_3<a_N-a_2$, the monomial $\frac{x^{a_3-a_1} \gamma_N (a_N-a_2) x^{a_N-a_3}}{a_2-a_1}$ is the required one.

    The computation for $\lambda_2$ is similar.

\end{proof}

\begin{lemma} \label{firstinfty}
	Consider a Laurent polynomial $f \in \C^{A \setminus a_1}$. Suppose that $f$ has no multiple roots at points in the set $\C \setminus 0$ and $f \notin \C^{A \setminus a_n}$.
	
	Then $f\in D_A$ only if $f \in \C^{A\setminus \{a_1,a_2\}}$.
\end{lemma}

\begin{proof}
	The Laurent polynomial $f$ has coefficients $(\gamma_1=0,\gamma_2,\gamma_3,\dots,\gamma_n)$. Suppose that $f\in D_A$. 
    
    Consider a $2$-dimesional affine subspace $W$ of $\C^A$ passing through $f$ transversally to $\C^{A\setminus\{a_1,a_2\}}$. It can be parameterized by \eqref{subspaceparam0} for some $\alpha$, $\beta$.
	
	For generic $W$ (such that the projection of $W$ onto $\C^{\{a_1,a_2\}}$ is surjective) one may further assume that 
    \begin{equation*} 
    \alpha_1=1, \beta_1 = 0,\alpha_2=0, \beta_2=1. 
    \end{equation*}
    
    Note that $(\lambda_1,\lambda_2)$  is a coordinate system on $W$. We obtain the  parameterization \eqref{sectionparam1} of the curve $D_A \cap W$. Similarly to part 2 of Lemma \ref{rationalsections} we obtain

	\begin{equation*} 
		\begin{split}
			\lambda_1(x) &= \frac{(a_3 \gamma_3 + \beta_3 a_2 \gamma_2 - a_2  \gamma_3 - a_3 \beta_3 \gamma_2) x^{a_2+a_3} +  o(x^{a_2+a_3})}{(a_2-a_1)x^{a_1+a_2} + o(x^{a_1+a_2})} =\\ &=  \frac{(a_3-a_2)(\gamma_3-\gamma_2 \beta_3)}{a_2-a_1}x^{a_3-a_1} + o(x^{a_3-a_1}) \\
			\lambda_2(x) &= \frac{(- a_2 + a_1 )\gamma_2 x^{a_1+a_2} + (- a_3 +a_1)\gamma_3x^{a_1+a_3} + o(x^{a_1+a_3})}{(a_2-a_1)x^{a_1+a_2} + o(x^{a_1+a_2})} =\\ &= -\gamma_2 - \frac{(a_3-a_1)\gamma_3}{a_2-a_1}x^{a_3-a_2} + o(x^{a_3-a_2}) 
		\end{split}
	\end{equation*}

Since $\gcd(A-A)=1$, for generic $W$, a generic element of $W \cap D_A$ has a unique multiple root (which is $x$). This follows from Proposition \ref{mainrank} (1). Thus, the parameterization is generically injective.

Now, since $f$ has no multiple roots and $f\notin \C^{A\setminus a_n}$, Lemma \ref{neighborhoods} implies that the unique branch of $D_A \cap W$ passing through $f$ corresponds to $x\to 0$. 

Hence, $\gamma_2 = \lim_{x \to 0} \alpha_2 \lambda_1(x) + \beta_2 \lambda_2(x) + \gamma_2 = 0 - \gamma_2 + \gamma_2 = 0$ and thus $f\in \C^{A\setminus \{a_1,a_2\}}$.
\end{proof}

\begin{defin}
    Denote the sets $\tilde B_1 = \{a_i-a_1 \mid i \ge 3\}$, $\tilde B_2 = \{a_i-a_2\mid i \ge 3\}$. For each nonnegative integer $r$, define the positive integer \[j_r=\gcd \bigcup_{i=1,2} \left(\tilde B_i \cap [a_3-a_i,a_3-a_i+r]\right).\]

    If $A$ is reduced, then $j_r=1$ for $r\ge a_n-a_1$.
\end{defin}

\begin{theor} \label{transversaltypeinfty}
	Consider a general point $f\in \C^{A \setminus \{a_1,a_2\}}$. Let $W$ be a general 2-dimensional affine subspace of $\C^A$ that intersects $\C^{A \setminus \{a_1,a_2\}}$ transversally at $f$. Consider the curve $W \cap D_A$. 
    
    (0) This curve passes through $f$ (in particular, $f\in D_A$). 
	
	(1) If $a_3-a_2=1$, then $f$ is a smooth point of $W \cap D_A$. 
	
	(2) If $a_2 = a_1+1$ and $a_3 = a_1 + 3$ then $f$ is an ordinary cusp of $W \cap D_A$.
	
	(3) Otherwise, the $\delta$-invariant of the singularity at $f$ equals \[\frac{(a_3-a_1-1)(a_3-a_2-1)}{2} + \sum_{r=0}^{a_n-a_1} \frac{j_r-1}{2}. \]
	
	(4) The singularity at $f$ is a 0-nondegenerate sparse curve singularity. 
\end{theor}

\begin{proof}
	Similarly to Proposition \ref{firstinfty}, denote the coefficients of $f$ by $(\gamma_1=0,\gamma_2=0,\gamma_3,\dots,\gamma_n)$. Again, consider a 2-dimensional subspace $W$ given by \eqref{subspaceparam0} for generic $\alpha$, $\beta$ and assume that the projection of $W$ onto $\C^{\{a_1,a_2\}}$ is surjective. Thus, we may impose the condition $\alpha_1=\beta_2=1, \alpha_2=\beta_1=0$. Then Lemma \ref{rationalsections} gives the parameterization \eqref{sectionparam1} of the curve $W\cap D_A$.
	
	A generic polynomial $f\in \C^{A \setminus \{a_1,a_2\}}$ has exactly one multiple root (which is 0). Also, $\gamma_3 \ne 0$. Thus, similarly to Proposition \ref{firstinfty} the curve $W \cap D_A$ has a unique branch at the point $(0,0)\in W$ that corresponds to $x\to 0$.  (Note that $(0,0)$ are coordinates of $f$). It is easy to see from the local parameterization \eqref{sectionparam2} of $W\cap D_A$ that this branch is smooth if $a_3-a_2 = 1$ and is an ordinary cusp if $a_3 - a_1 = 3$ and $a_3 - a_2 = 2$.
	
	By \cite{esv}, the $\delta$-invariant of the singularity at $f$ is at least $(a_3-a_1-1)(a_3-a_2-1)/2$. We check the 0-nondegeneracy condition of \cite{esv} for $\lambda_1,\lambda_2$ to obtain the exact value of the $\delta$-invariant.

     Let $\kappa>1$ be a common divisor of $a_3-a_1$ and $a_3-a_2$. By $N$ denote the smallest integer such that $a_N-a_3$ is not divisible by $\kappa$.

    We have already computed the necessary coefficients of $\lambda_1$ and $\lambda_2$ in parts 2 and 3 of Lemma \ref{rationalsections}. The 0-nondegeneracy condition reads as follows:
    \begin{equation*}
        \frac{\gamma_N(a_N-a_2)}{a_2-a_1} \cdot \left( \frac{\gamma_3(a_3-a_2)}{a_2-a_1} \right)^{-1} \cdot (a_3-a_1)^{-1} \ne \frac{\gamma_N(a_1-a_N)}{a_2-a_1} \cdot \left( \frac{\gamma_3(a_1-a_3)}{a_2-a_1} \right)^{-1} \cdot (a_3-a_2)^{-1}.
    \end{equation*}

    Assume $\gamma_i\ne 0$ for $i\ge 3$. Then this inequality is reduced to $a_N-a_2\ne a_N-a_1$ which is trivial. 
\end{proof}

\begin{theor} \
    The variety $D_A \cap \C^{A\setminus a_1}$ has at most 2 irreducible components: $\C^{A\setminus\{a_1,a_2\}}$ and $D_{A\setminus a_1}$. 
    
    The first component belongs to the singular locus of $D_A$ if and only if one of the conditions (2), (3) of Theorem \ref{transversaltypeinfty} hold. 
    
    The second component is not contained in to $D_A^{sing}$ if $A$ is doubly reduced.
\end{theor}

\begin{proof}
    If a point $f \in D_A \cap \C^{A\setminus a_1}$ does not belong to $\C^{A\setminus\{a_1,a_2\}}$ then Lemma \ref{firstinfty} implies that $f$ has a singularity in $\C^*$ and thus belongs to $D_{A\setminus a_1}$. This variety is irreducible and thus the first statement of the theorem follows.

    The second statement is contained in Theorem \ref{transversaltypeinfty}.

    Now we proceed to the last statement. The set $A \setminus a_1$ is reduced. Thus, a general element $f$ of $D_{A\setminus a_1}$ has exactly one singularity which we denote by $p\in \C^*$. Thus, $f$ is a smooth point of $D_A$ by the same argument as in Proposition \ref{mainnotinfty}: apply Proposition \ref{stratadimension} (B) and the rank estimation of Theorem \ref{mainrank} (1) in the case $m=1$ and point $p$.
\end{proof}

\section{Singularities of surface projections}

We classify and enumerate singularities of the  projection $p:\{f=0\}\to\CC^{n-1}$ forgetting the last coordinate for a general Laurent polynomial $f:\CC^n\to\C$ with a given Newton polytope $N\subset\Z^n$.
The problem reduces to the study of the sparse discriminant $D_A\subset\C^A$ with $A:=h(N)$, via the map $\tilde f:\CC^{n-1}\to\C^A$ sending $x$ to the fiber $f(x,\cdot)$. Recall that $h:\Z^n\to\Z$ is the last coordinate, referred to as the {\it height}.

The first subsections employ the map $\tilde f$ to prove the results from the introduction, and to give a combinatorial formula for the number of stable multisingularities (double folds) in the case $n=3$. 
However, the latter formula is not practical: similarly to the classical singularity theory, the study of multisingularities in our setting is drastically harder than (mono)singularities. 

This motivates the rest of this section, counting double folds more practically under a mild additional assumption on the Newton polytope $N$ (normality). It is based on computing the ultratropicalization \cite{esv} of the discriminant of $p$. 

\subsection{Types of singularities} 

\begin{defin}
For a smooth map $F:M\to N$, the {\it Thom--Boardman stratum} $S_1\subset M$ is the critical locus, and $S_i$ is the rank drop locus of $dF$ restricted to $S_{i-1}$ (provided that the latter is smooth).
\end{defin}
\begin{exa}\label{exatb}
1. For a stable map of surfaces, $S_2$ is the set of cusps.

2. (see Proposition \ref{whitneycusp}). For the unfolding $\{(s,g),\,G(s):=s^k+g(s)=0\}\to\C[s]/s^k,\,(s,g)\mapsto [g]$, the Thom--Boardman stratum $S_i$ is the set of pairs $(s,g)$ such that $s$ as a root of $G$ has multiplicity at least $i+1$, i.e.
$$S_i=\{(s,g)\,|\,G(s)=G'(s)=\cdots=G^{(i)}(s)=0\}.$$
\end{exa}
\begin{theor}\label{thtb}
In the setting at the beginning of this section, the $i$-th Thom--Boardman stratum of $p$, for $i< |h(N)|$ and general $f\in\C^N$, is given by the equations
$$f(x)=\partial f/\partial x_n=\cdots=\partial^i f/\partial x_n^i=0.\eqno{(*)}$$
\end{theor}
\begin{proof}
Let $x$ satisfy $(*)$ and $\partial^{i+1} f/\partial x_n^{i+1}\ne 0$. Applying Lemma \ref{deformrankV} 
to $f_0=\tilde f(p(x_0))$ and $x_0:=x_n$, we get a full rank analytic map $L:(\C^{h(N)},f_0)\to(\C[s]/s^i,0)$, so that the subvariety $V:=L^{-1}(0)$ is smooth near $f_0$. By Corollary 6.7 in \cite{esv}, the map $\tilde f$ is transversal to $V$ at $p(x_0)$, so the composition $L\circ\tilde f$ has the full rank, so the Thom--Boardman strata of $p$ near $x_0$ induce from the ones from Example \ref{exatb}.2.
\end{proof}

More generally, Corollary 6.7 in \cite{esv} implies that the map $\tilde f$ for general $f$ is transversal to the Whitney strata of the discriminant $D_A$. Applying this to the strata $U$, $D_\prec$, $D_\times$, described in Theorem \ref{mainth}, we get the following.

\begin{sledst}\label{csingp} Assume $n=3$. 

1. If $|h(N)|>2$, the ordinary cusps and folds are the only singularities of $p$.

2. If moreover $h(N)$ satisfies Convention \ref{conven1}, then double folds are the only further multisingularities of $p$.
\end{sledst}
\begin{rem}
For now, we do not analyze the singularities of the diescriminant of $p$ at the non-properness points of $p|_{\mbox{\scriptsize critical locus of }p}$. 
\end{rem}
\begin{proof}
In the notation of Theorem \ref{mainth}, the map $\tilde f$ is transversal to $U$, $D_\prec$, $D_\times$ and $D_A\cap\{c_a=0\}$ for every $a\in A$. Since $p$ is the pull back of $p_A$ under $\tilde f$, the multisingularities of $p$ over the preimages of the listed strata are the same as the respective multisingularities of $p_A$ (described in Theorem \ref{mainth} over $U$, $D_\prec$ and $D_\times$). 
Since the rest of $D_A$ (over which $p_A$ has worse singularities) has codimension greater than 2 in $\C^A$, its preimage under $\tilde f$ is empty.
\end{proof}

\subsection{Tropical enumeration of singularities} 

We now explain how to count these (multi)singularities. The coefficients of the generic polynomial $\sum_{a\in A}c_at^a$ form a coordinate system on $\C^A$, and we accordingly decompose $f(x,t)=\sum_{a\in A}f_a(x)t^a$. The graph $\Gamma_f$ of the map $\tilde f:\CC^2\to\C^A$ is the complete intersection of the hypersurfaces $c_a-f_a(x)=0$, whose Newton polytopes we denote by $N_a$ (combinatorially, this polytope spans the point corresponding to the term $c_a$ and a shifted copy of the set $h^{-1}(a)\cap N$).

For general $f$, the intersection number of $\Gamma_f$ with any given semialgebraic set $H\subset\CC^2\times\CC^A$ of complimentary codimension (2, in our case) equals the tropical intersection number of the dual tropical fans $[N_a]$ of the Newton poyltopes and  $\Trop H$ (the tropicalization of the closure of $H$). 
We shall apply this to the strata $D_\prec$ and $D_\times$ (described in Theorem \ref{mainth}) to count the singularities of the projection $p$. 

\begin{theor}
For $n=3$, assume $h(N)$ satisfies convention \ref{conven1}. Then the number of cusps and double folds of the projection $p$ for general $f$ equals the tropical intersection numbers $(\R^2\times\Trop D_\prec)\cdot\prod_{a\in A}[N_a]$ and $(\R^2\times\Trop D_\times)\cdot\prod_{a\in A}[N_a]$ in $\R^2\times\R^A$, respectively.
\end{theor}
This is a combinatorial expression for the number of singularities, because the tropical fans $\Trop D_\prec$ and $\Trop D_\times$ are computed in terms of $A$ in \cite{jems} and \cite{alicia}. 
\begin{proof}
For general $f$, the graph $\Gamma_f$ is a general complete intersection with the Newton polytopes $N_a$. (Rigorously: if a property of a complete intersection is invariant under multiplying its equations by non-zero numbers, then it is satisfied for the complete intersection of general hypersurfaces with the Newton polytopes $N_a$ iff it is satisfied for $\Gamma_f$ corresponding to general $f\in\C^N$.)

Thus the intersection number of any codimension two set $H\subset\CC^2\times\CC^A$ with $\Gamma_f$ equals its intersection number  with a general complete intersection whose Newton polytopes are $N_a$, i.e. $(\Trop H)\cdot\prod_{a\in A}[N_a]$.

Setting $H$ to $D_\prec$ and $D_\times$ we count all cusps and double folds of $p$ at the points of $\CC^2$ sent to $\CC^A$ by $\tilde f$. The points of $\CC^2$ missing $\CC^A$ cannot be multisingularities of $p$ under our assumption on $A$, by Theorem \ref{mainth}.2.
\end{proof}

In practice, computing the tropical fans of the singularity strata is difficult, so we give alternative formulas for the number of singularities, resembling the Bernstein--Kouchnirenko--Khovanskii formula.

\subsection{Enumeration of cusps: sch\"on complete intersections.} 

By Theorem \ref{thtb}, a Thom--Boardman stratum of the projection $p$ is the complete intersection $f=\partial f/\partial x_n=\partial^2 f/\partial x_n^2=\cdots=0$ for a general polynomial $f(x)=\sum_{b\in B} c_bx^b$ with the support set $B\subset\Z^n$. We employ the fact that it belongs to a certain convenient class introduced in \cite{eci}.

A complete intersection of hypersurfaces $f_i=0$ in the torus $\CC^n$ is {\it nondegenerate upon cancellations} ({\it NUC}, or, equivalently, is {\it sch\"on}  and {\it Newtonian}), if every linear function $\gamma:\Z^n\to\Z$ admits {\it cancellations} $\tilde f_i:=f_i+c_{i,i-1}f_{i-1}+\cdots+c_{i,1}f_1$ (with coefficients $c_{i,j}$ Laurent polynomials) whose $\gamma$-leading parts form a regular complete intersection. Denoting the $\gamma$-degree of $\tilde f_i$ by $m_i(\gamma)$, we get conewise linear functions $m_i$ (each defined on the corner locus of the preceding one), which are called {\it tropicalizations} of the initial equations $f_i$, and in what follows are always arbitrarily extended to conewise linear continuous functions on the whole $\R^n$. 

Geometry of a NUC complete intersection can be expressed in terms of the tropicalizations $m_i$. Writing these expressions is simplified with the following formalism.
\begin{defin}
Given a homogeneous degree $n$ polynomial $P\in\R[\mu_1,\ldots,\mu_k]$ and conewise linear functions $m_1,\ldots,m_k$ on $\R^n$, the value $P(m_1,\ldots,m_k)\in\R$ is defined by the following properties:

1. For fixed $(m_1,\ldots,m_k)$, the map $P\mapsto P(m_1,\ldots,m_k)$ is linear in $P$.

2. If $Q(\mu_1,\ldots,\mu_k)=P((\mu_1,\ldots,\mu_k)\cdot C)$ for a $k\times k$-matrix $C$, then $Q(m_1,\ldots,m_k)=P((m_1,\ldots,m_k)\cdot C)$.

3. The $n$-th power of a convex conewise linear function $m$ equals the lattice volume of the polytope $N_m$ with the support function $m$ (i.e. the Newton polytope of $m$ regarded as a tropical polynomial).

More generally, for any polynomial $P$, its value $P(m_1,\ldots,m_k)$ in the ring of tropical fans in $\R^n$ is defined by the same properties (1) and (2) and

3'. The $j$-th power of a convex conewise linear function $m$ equals the $(n-j)$-dimensional skeletone of the dual fan of the polytope $N_m$, in which the weight of every maximal cone is defined as the lattice volume of the dual $j$-dimensional face of $N_m$.
\end{defin}
\begin{rem}\label{remmixed}
1. The first definition embeds in the second one by identifying a number $c\in\R$ with the 0-dimensional tropical fan $\{0\}$ bearing weight $c$.

2. In this formalism, the product of several convex functions $m_i$ equals the intersection product of their dual fans, and, in particular, the product of $n$ such functions $m_i$ equals the mixed volume of their polytopes $N_{m_i}$.

3. We can "naively" plug the conewise linear functions $m_i$ in to the polynomial $P$ pointwise, to get a conewise polynomial function $\varphi:\R^n\to\R$. By \cite{??mmj}, the value of $P(m_1,\ldots,m_k)$ can be read off from $\varphi$ (by iteratevely applying the corner locus operator).
\end{rem}
\begin{theor}\label{thnuc}(\cite{eci})\label{nucth}
The Euler characteristics and the tropical fan of the set $f_1=\cdots=f_k=0$ are, respectively, the coefficients of $t^n$ and $t^k$ in the power series expansion of 
$$\frac{tm_1}{1+tm_1}\cdot\ldots\cdot\frac{tm_k}{1+tm_k}.$$
\end{theor}
\begin{rem}
The other coefficients in this series are so called {\it tropical characteristic classes} of $f=0$ \cite{jems}, \cite{gb}.
\end{rem}

For many NUC complete intersections $f$ (such as engineered and simplest symmetric ones), the sum $m_1+\cdots+m_i$ is convex, i.e. a support function of a certain lattice polytope \cite{eci}. This polytope is called the $(i-1)$-th {\it incremental polytope} $N^{(n-1)}$ of $f$. 

For the Thom--Boardman complete intersection of Theorem \ref{thtb}, the $i$-th incremental polytope $N^{(i)}$ can be represented (by \cite{eci} and \cite{kks} respectively) as

$$\bigcap_{q_1,\ldots,q_i\in\Z} \conv\Bigl( N+(N\setminus h^{-1}q_1)+\cdots+(N\setminus h^{-1}\{q_1,\ldots,q_i\})\Bigr)=$$
$$=\conv\{a_0+\cdots+a_i\,|\,a_0,\ldots,a_i\in N,\,h(a_0),\ldots,h(a_i)\mbox{ are pairwise different}\}.$$
In particular, $N^{(0)}=N$ and $N^{(1)}=N'$. Plugging the differences of the support functions these consecutive incremental polytopes as $m_i$'s in Theorem \ref{thnuc}, and evalutating in terms of their mixed volumes and dual fans by Remark \ref{remmixed}, we get the following.
\begin{theor}\label{thtbenum}
In the setting of Theorem \ref{thtb}, the Euler characteristics and the tropical fan of the Thom--Boardman stratum $(*)$ are, respectively, the coefficients of $t^n$ and $t^{i+1}$ in the power series expansion of 
$$\frac{tN}{1+tN}\cdot\frac{t(N'-N)}{1+t(N'-N)}\cdot\ldots\cdot\frac{t(N^{(i)}-N^{(i-1)})}{1+t(N^{(i)}-N^{(i-1)})}.$$
\end{theor}

\subsection{Enumeration of nodes: ultratropicalizations.} 
\begin{assum}
A spanning polytope $N\subset\Z^3$ is {\it normal}, if the projective $N$-toric variety is normal. Combinatorially, this means that, for each face $F\subset N$, the affine surjection $\Z^3\to\Z^{\codim F}$, sending $F$ to 0, sends $N$ to a generating set of a saturated semigroup $S$ (i.e. the set of all lattice points in a rational polyhedral cone).

The union of bounded faces of $\conv(S\setminus\{0\})$ is called the {\it link diagram} $N_F$ of $N$ at $F$.
\end{assum}
\begin{rem}\label{remdiagr}
1. Non-normal polytopes in $\Z^3$ do exist (e.g. a tetrahedron with vertices $(000),(001),(111)+k(221),(-110)-k(-221)$ for $k>0$ at the edge $F$ spanning the first two vertices), but are so rare that probably admit a tractable classification.

2. If $F$ is an edge of normal $N$, then each edge of its link diagram $N_F$ contains all lattice points of its convex hull and has unit lattice distance from 0.

If $F$ is a nonhorizontal 2-face of normal $N$, then the only point of its link diagram has unit lattice distance from 0, and the preimage $F'\subset N$ of this point has the following properties:

-- each of $h(F)$ and $h(F')\subset\Z$ can be shifted to $d_F\Z$, but $h(F)$ cannot be shifted to its proper sublattice;

-- $h(F\cup F')$ cannot be shifted to a proper sublattice of $\Z$.
\end{rem}

To express the number of double folds of the projection $p$ in terms of the polytope $N$, introduce the following notation.

1. If $F\subset N$ is a {\it horizontal} edge (such that $h(F)$ is a point), write $h_F$ for the unique element of $h(F)$ and define its {\it numerator} $n_F$ as its lattice length.

2. For a horizontal edge $F$, let $\pi_F$ be the affine lattice projection contracting $F$ to 0. Define $N'_F$ as the Newton diagram of the following convex hull, i.e. the part of its boundary visible from 0:
$$
\conv\pi_F\bigl(N\setminus h^{-1}(h_F)\bigr).
$$
Thus the height-$h_F$ layer is deleted before taking the convex hull. We shift the height induced on the quotient lattice by $h_F$. For every edge $E\subset N'_F$ on which the shifted height is nonconstant and whose affine span has positive lattice distance from 0, define
$$
n_E:=\mathop{\rm length}\nolimits_{\Z}h(E),\qquad
d_E:=\frac{n_E}{\mathop{\rm length}\nolimits_{\Z}E},\qquad
\rho_E:=\dist_{\Z}(0,\mathop{\rm aff}E).
$$
Edges on which $h$ is constant will not be used. An exterior normal covector of the flag $(F,E)$ projects to a nonzero ray $R$ in the two-dimensional base. We denote by $\mathcal E_R(F)$ the set of the above edges with the same projected exterior-normal ray $R$.

Every edge of $N'_F$ contained in an edge of $N_F$ has $\rho_E=1$. The deletion may additionally create at most one edge $E_0$ crossing the height-zero line. If it does, then its end points have heights $-1$ and $1$, so $n_{E_0}=2$ and $d_{E_0}\in\{1,2\}$, while $\rho_{E_0}$ may be arbitrary. Some edges of $N_F$ may also be shortened or disappear upon the deletion.

3. If $F\subset N$ is a nonhorizontal 2-face with a primitive normal vector $(\gamma_1,\gamma_2,\gamma_3)$, then $n_F$ is the lattice area of $F$ minus the lattice lengths of its horizontal edges (if any), and
$$
d_F=\gcd(\gamma_1,\gamma_2)=\gcd\{h(a)-h(b)\mid a,b\in F\}.
$$

4. Define the numbers $h_\pm\in\Z$ as the end points of the segment $[h_-,h_+]=h(\conv N)$, and $V_\pm\in\Z$ as the lattice mixed area of the convex hulls of the polygons $N\cap h^{-1}(h_\pm)$ and $N\cap h^{-1}(h_\pm\mp1)$ (the two top or respectively bottom horizontal layers of $N\subset\Z^3$). We put $V_\pm=0$ if the respective adjacent layer is empty.

5. Assuming as always $|h(N)|>2$, put
$$
r_+=h_+-\max\bigl(h(N)\setminus\{h_+,h_+-1\}\bigr),
$$
$$
r_-=\min\bigl(h(N)\setminus\{h_-,h_-+1\}\bigr)-h_-.
$$
The sets in these formulas are nonempty by the standing assumption $|h(N)|>2$.

We denote the logarithmic derivative $x_i\partial f/\partial x_i$ by $\partial_i f$. Let $C$ be the critical locus of $p$, let $D\subset\CC^2$ be the reduced closure of $p(C)$, let $K$ be the number of cusps of $p$, and let $\Delta_D$ be the Newton polygon of $D$.

\begin{theor}\label{enumnodes}
If all singularities of the projection $p$ are stable (which happens e.g. under the assumptions of Corollary \ref{csingp}), and $N$ is normal, then the number of double folds of $p$ equals
$$
\begin{aligned}
\frac12\Bigg[&-\sum_{\epsilon=\pm}(r_\epsilon^2-3r_\epsilon+1)V_\epsilon
-\sum_{\substack{F\ {\rm a\ nonhorizontal}\\ {\rm two\mbox{-}face}}}n_F(d_F-1)\\
&-\sum_{\substack{F\ {\rm a\ horizontal\ edge}\\R\ {\rm a\ base\ ray}}}
\sum_{E\in\mathcal E_R(F)}
n_Fn_E(n_E-1)\frac{\rho_E}{d_E}\\
&-2\sum_{\substack{F\ {\rm a\ horizontal\ edge}\\R\ {\rm a\ base\ ray}}}
\sum_{\substack{\{E_1,E_2\}\subset\mathcal E_R(F)\\E_1\ne E_2}}
n_Fn_{E_1}n_{E_2}
\min\left(\frac{\rho_{E_1}}{d_{E_1}},\frac{\rho_{E_2}}{d_{E_2}}\right)\\
&+\chi(C)-2K+\mathop{\rm Area}\nolimits_{\Z}(\Delta_D)\Bigg].
\end{aligned}
$$
The innermost sum is over unordered pairs of distinct edges.
\end{theor}

\begin{proof}
Let $M$ be the number of double folds and put ${P}:=D\setminus p(C)$.
We first determine the local invariants of the points of ${P}$. These points are projections of the points of the closure of $C$ in $\CC^2\times\CP^1\supset\CC^3$ that belong to $\CC^2\times\{0\}$ and $\CC^2\times\{\infty\}$ respectively. We refer to these two kinds of points as the {\it bottom} and {\it top} points, because their geometry is governed by the geometry of the polytope $N$ near the planes $h^{-1}(h_-)$ and $h^{-1}(h_+)$ in the sense made precise below. We study top points, the case of bottom points is identical. 

Let $q$ be the least positive integer such that
$N\cap h^{-1}(h_+-q)\ne\varnothing$.  If the mixed area of the top layer of $N$ and
the height-$(h_+-q)$ layer is nonzero, then $q=1$.  Indeed, this follows from
saturation in the rank-one quotient when the top face is two-dimensional.  If
the top face is one-dimensional, nonzero mixed area gives two consecutive
points $(a,q),(a+1,q)$ in the first nonzero slice of the rank-two quotient.
For $q>1$, their cone contains a lattice point of height strictly between $0$
and $q$; saturation puts this point in the semigroup generated by the image of
$N$, contradicting the minimality of $q$.  If the top face is a point, the
mixed area is zero.

The equations defining the closure of the critical curve show that every top
point of ${P}$ projects to a common zero of the restrictions of $f$ to the
top layer and the first occupied lower layer.  Thus, by the BKK theorem,
the top points of ${P}$ are precisely the $V_+$ transversal common zeros
of the restrictions $f_0,f_1$ to the height-$h_+$ and height-$(h_+-1)$ layers.
At each such zero, the restriction $f_r$ to the height-$(h_+-r_+)$ layer does
not vanish for general $f$.

Fix one of these zeros, put $r=r_+$, take $s=x_3^{-1}$, and multiply $f$ by
$s^{h_+}$.  Since $u=f_1$ and $v=f_0$ are local coordinates in the base, the
equation has the form
$$
v+su+s^r(c+\psi(u,v,s))=0,
\qquad c\ne0,\qquad \psi(0,0,0)=0.
$$
The closure of the critical curve is defined by this equation and its
$s$-derivative.  Their Jacobian with respect to $(u,v)$ is nonsingular at the
origin; hence the implicit function theorem gives its unique parametrization
$$
u=-rcs^{r-1}+O(s^r),\qquad
v=(r-1)cs^r+O(s^{r+1}).
$$
Since $\gcd(r-1,r)=1$, its image is an irreducible plane branch with equation
$$
u^r+\lambda v^{r-1}+\mbox{terms of strictly greater $(r-1,r)$-weight}=0,
\qquad \lambda\ne0.
$$
Its Milnor number is $(r_+-1)(r_+-2)$.  The identical argument at the bottom
gives $V_-$ branches with Milnor number $(r_--1)(r_--2)$.  Consequently
$$
|{P}|=V_++V_-,\qquad
\sum_{a\in{P}}(1-\mu_a)
=-\sum_{\epsilon=\pm}(r_\epsilon^2-3r_\epsilon+1)V_\epsilon.
$$

No point of ${P}$ belongs to $p(C)$.  At the top, such a coincidence
would give a solution of the ECI given by four equations and three variables
$$
f_0=f_1=f=\partial_3f=0,
$$
which is inconsistent for general $f$ by \cite{eci}.
A top point cannot coincide with the bottom point by a similar dimension count.

Let $D^\nu$ be the normalization of $D$.  Stability identifies $C$ with the
complement in $D^\nu$ of the normalization points over ${P}$.  The local
calculation and the preceding noncollision statements give one such point over
each $a\in{P}$, and hence
$$
\chi(D^\nu)=\chi(C)+|{P}|.
$$
All other singularities of $D$ are the $M$ nodes corresponding to double folds
and the $K$ ordinary cusps.  Since every germ at a point of ${P}$ is
unibranch,
$$
2\sum_{z\in\mathop{\rm Sing}D}\delta(D,z)
=2M+2K+\sum_{a\in{P}}\mu_a.
$$

For every ray $R$ of $\mathop{\rm Trop}D$, let $(g^R_{ij})$ be its tangency
matrix, and put $G:=\sum_{R,i,j}g^R_{ij}$.  Lemma 1.11 in \cite{esv} and the
preceding identities give
$$
2M=\chi(C)+\mathop{\rm Area}\nolimits_{\Z}(\Delta_D)-2K-G
-\sum_{\epsilon=\pm}(r_\epsilon^2-3r_\epsilon+1)V_\epsilon.
$$

In the $d$-fold pullback of Theorem \ref{thut}, an intersection number $m$
contributes $m/d$ to the corresponding off-diagonal entry of the tangency
matrix, while its diagonal entries are zero.  Counting ordered pairs of the
branches listed in that theorem gives
$$
\begin{aligned}
G={}&\sum_{\substack{F\ {\rm a\ nonhorizontal}\\ {\rm two\mbox{-}face}}}
n_F(d_F-1)\\
&+\sum_{\substack{F\ {\rm a\ horizontal\ edge}\\
R:\,{E}_R(F)\ne\varnothing}}
\sum_{E\in{E}_R(F)}
n_Fn_E(n_E-1)\frac{\rho_E}{d_E}\\
&+2\sum_{\substack{F\ {\rm a\ horizontal\ edge}\\
R:\,{E}_R(F)\ne\varnothing}}
\sum_{\substack{\{E_1,E_2\}\subset{E}_R(F)\\E_1\ne E_2}}
n_Fn_{E_1}n_{E_2}
\min\left(\frac{\rho_{E_1}}{d_{E_1}},
          \frac{\rho_{E_2}}{d_{E_2}}\right).
\end{aligned}
$$
The innermost sum is over unordered pairs.  Substitution of this expression
for $G$ into the preceding formula and division by $2$ gives the stated formula.
\end{proof}

\section{The ultratropicalization of the discriminant}
This is the only ingredient in the proof of Theorem \ref{enumnodes} that have not been settled in the previous sections. In the setting of this theorem, fix a ray $R$ in the character lattice $\Z^2$ of the base torus $\CC^2$ and make a monomial change of the two base coordinates which identifies it with $(-1,0)$. In this toric chart, let $D\subset\C\times\C^*$ denote the reduced closure of the discriminant curve, and let $\tilde D$ be its reduced preimage under the endomorphism $(x_1,x_2)\mapsto(x_1^d,x_2)$. For sufficiently divisible $d$, all local components $\tilde D_i$ of $\tilde D$ near $L=\{x_1=0\}$ are smooth and transversal to it.

Describing the ultratropicalization of $D$ in the direction $R$ amounts to enumerating the components $\tilde D_i$ and computing their pairwise orders of tangency (or, equivalently, intersection numbers).
\begin{theor}\label{thut} In the setting of Theorem \ref{enumnodes}, assume $N$ is normal and $|h(N)|>2$.
Choose a sufficiently divisible integer $d>0$. Then the local components $\tilde D_i$ along $L$ can be indexed by the following tuples:

1) $(F,i,j)$, where $F$ is a nonhorizontal 2-face whose exterior-normal ray projects to $R$, $i\le n_F$, and $j\le d_F$;

2) $(F,i,E,j)$, where $F$ is a horizontal edge, $E\in\mathcal E_R(F)$, $i\le n_F$, and $j\le n_E$;

in such a way that the local intersection numbers of distinct components along $L$ are as follows:

-- components $(F,i,j)$ and $(F,i,j')$, where $j\ne j'$, have intersection number $d/d_F$;

-- components $(F,i,E,j)$ and $(F,i,E',j')$, where $(E,j)\ne(E',j')$, have intersection number
$$
d\min\left(\frac{\rho_E}{d_E},\frac{\rho_{E'}}{d_{E'}}\right);
$$

-- no other two components meet along $L$.
\end{theor}

To construct this indexing, set
$$
f_k(x_1,x_2,x_3):=f(x_1^d,x_2,x_3x_1^k),
$$
and let $C_k$ be the closure 
of the curve $f_k=\partial_3f_k=0$ in $\C\times\CC^2$ (equivalently, its ideal is saturated by $x_1$). For $d$ divisible enough and general $f$, all local components of $C_k$ near the plane $H=\{x_1=0\}$ are smooth, and intersect it transversally at pairwise different points $Y_k=C_k\cap H$.

The components $C_{k,y},\,y\in Y_k$, project to the sought components of $\tilde D$, because it is the discriminant of the projection $p:\{f_k=0\}\to\CC^2$ regardless of $k$. Conversely, after making $d$ divisible by all Puiseux denominators, every local branch of the normalization of the critical curve has integral order $k$ in the third coordinate. Rescaling this coordinate gives a branch of some $C_k$ meeting $H$. Its initial equations are supported either on a nonhorizontal 2-face $F$, or on a horizontal edge $F$ followed, after deletion of its height layer, by an edge $E\subset N'_F$. Thus the above list is exhaustive. Transversality gives a unique germ through every point of $Y_k$, and Lemma \ref{mainl} below computes its projection and excludes duplicate labels.

So the proof of Theorem \ref{thut} amounts to enumerating the components $C_{k,y}$ and computing the orders of tangency between their projections. The answer is summarized below.

\begin{lemma}\label{mainl}
1. If the covector $(-d,0,-k)$ is exterior normal to a nonhorizontal 2-face $F\subset N$, then the set $p(Y_k)$ has $n_F$ points, and $Y_k$ has $d_F$ points
$$
\{(x_2,\zeta_{d_F}^jx_3)\mid j=0,\ldots,d_F-1\}
$$
over each $x_2\in p(Y_k)$, where $\zeta_{d_F}$ is a primitive $d_F$-th root of unity.

2. If the covector $(-d,0,-k)$ exposes a horizontal edge $F\subset N$ and an edge $E\in\mathcal E_R(F)$ of $N'_F$, then the set $p(Y_k)$ has $n_F$ points, and $Y_k$ has $n_E$ points over each of them.

3. The sets $p(Y_k)$ do not overlap for different $k\ne k'$, unless both are as in (2) with the same horizontal edge $F$, but different edges $E$ and $E'$ of $N'_F$. In the latter case, $p(Y_k)=p(Y_{k'})$.

4. Assume $k$ is as in (1). For distinct $y,y'\in Y_k$, the components $p(C_{k,y})$ and $p(C_{k,y'})$ have intersection number $d/d_F$ if $p(y)=p(y')$ and 0 otherwise.

5. Assume that $C_{k,y}$ and $C_{k',y'}$ are distinct components as in (2), corresponding to the same horizontal edge $F$ and to not necessarily distinct edges $E,E'\subset N'_F$. Their projections have intersection number
$$
d\min\left(\frac{\rho_E}{d_E},\frac{\rho_{E'}}{d_{E'}}\right)
$$
if $p(y)=p(y')$, and 0 otherwise.

6. Except for (4) and (5), distinct local components $p(C_{k,y})$ do not meet along $L$.
\end{lemma}

The rest of the section is devoted to the proof of this lemma.

\subsection{Proof of Lemma \ref{mainl}.1-3,6: engineered complete intersections.}

Given a (multi)set $A\subset\Z^n$ ({\it support set}) and a rank $k$ matrix $V=(v_{a,i})$ ({\it kinetic matrix}) an {\it engineered complete intersection} is a system of equations $f_i(x):=\sum_{a\in A}v_{a,i}c_ax^a=0,\,i=1,\ldots,k$, where the vector of coefficients $(c_a,\,a\in A)$ is generic.
\begin{rem}
A statement $P$ about an ECI with a given support set and kinetic matrix should be interpreted as a statement that there exists a Zariski open set $U\subset\C^A$ such that $P$ holds for all coefficient vectors $(c_a)\in U$.
\end{rem}
\begin{exa}
If the Laurent polynomial $f$ is a generic linear combination of monomials from $A\subset\Z^n$, then any system whose equations are logarithmic derivatives of $f$ and/or are obtained from $f$ by deleting some monomials is an ECI, provided that the corresponding kinetic rows have full rank.
\end{exa}
\begin{theor}[\cite{eci}]\label{theci}
1. Every ECI is a regular NUC complete intersection in $\CC^n$.
In particular, if $A$ can be shifted to a less than $k$-dimensional plane, the ECI has no solutions in $\CC^n$.

2. The $k$-th incremental polytope of the ECI $f$ depends only on the \textbf{matroid complete intersection} $(A,M)$, where $M$ is the matroid on $A$ represented by the kinetic matrix $V$.
\end{theor}

\begin{defin}
Recall that, for a polynomial $f(x)=\sum_{a\in\Z^n}c_ax^a$,

-- the {\it support set} is the set of all $a$ such that $c_a\ne0$;

-- the {\it restriction} $f|_F$ onto a finite set of monomials $F\subset\Z^n$ is defined as $\sum_{a\in F}c_ax^a$.
\end{defin}

{\it Proof of Lemma \ref{mainl}.1.} i. Let $\tilde g(x_2,x_3)$ be the lowest coefficient of a monomial normalization of $f_k$ as a Laurent polynomial of $x_1$. From the setting, up to a monomial multiplier,
$$
\tilde g(x_2,x_3)=g(x_2,x_3^{d_F}).
$$
The support set $B$ of $g$ is lattice-isomorphic to $F\subset N$, and the respective monomial change of coordinates identifies $g$ with the restriction $f|_F$.

ii. From the definition of $d_F$ and Remark \ref{remdiagr}, $h(B)$ cannot be shifted to a proper sublattice of $\Z$.

iii. We want to prove that the system
$$
g=\partial_3g=0\eqno{(*)}
$$
has $n_F$ roots, and they have pairwise distinct values of coordinate $x_2$. This will imply the statement of the lemma, because $Y_k=\{\tilde g=\partial_3\tilde g=0\}$ is the preimage of $(*)$ under
$$
(x_2,x_3)\longmapsto(x_2,x_3^{d_F}).
$$

iv. The root count for $(*)$ is Theorem 1.8 in \cite{eci}, because $(*)$ is an ECI supported at $B\simeq F$. It remains to prove that the roots have pairwise different $x_2$. These values of $x_2$ are the preimages of the discriminant $D_{h(B)}\subset\C^{h(B)}$ under the map
$$
x_2\longmapsto g(x_2,\cdot).\eqno{(**)}
$$
For general $g=f|_F$, the map $(**)$ intersects the discriminant at its general points by Corollary 6.7 of \cite{esv}. Since $h(B)$ cannot be shifted to a proper sublattice, a general polynomial in this discriminant has only one multiple root. Thus every such $x_2$ contributes only one root of $(*)$. $\hfill\square$

\vspace{1ex}

{\it Proof of Lemma \ref{mainl}.2.}  Put $t=x_1$. After multiplying $f_k$ by a monomial, including $x_3^{-h_F}$, write it near the flag $(F,E)$ as
$$
\tilde f_k(t,x_2,x_3)=G_F(t,x_2)+t^{q_E}H_E(x_2,x_3)+O(t^{q_E+1}),
\qquad q_E=\frac{d\rho_E}{d_E}.
$$
Here $G_F$ is the contribution of the complete height-$h_F$ layer, so it is independent of $E$, and $g_F(x_2):=G_F(0,x_2)$ is $f|_F$ up to a monomial multiplier. The critical equation is unchanged by this normalization after replacing $\partial_3f_k$ by the elementary row operation $(\partial_3-h_F)f_k$. Since $C_k$ is the saturated closure, its limiting equations are
$$
g_F(x_2)=0,\qquad \partial_3H_E(x_2,x_3)=0.\eqno{(*)}
$$
Indeed,
$$
t^{-q_E}\partial_3\tilde f_k=\partial_3H_E+O(t).
$$

The first equation of $(*)$ is general of degree $n_F$, and hence has $n_F$ distinct roots in $\C^*$. The height-zero layer was deleted before $N'_F$ was formed, so logarithmic differentiation does not delete either end point of the $x_3$-support of $H_E$. Therefore the second equation has $n_E$ simple roots in $\C^*$ after every root of $g_F$ is substituted. This proves the statement, because the projection $p$ forgets the coordinate $x_3$. $\hfill\square$

\vspace{1ex}

{\it Proof of Lemma \ref{mainl}.3 (and 6).} We need to consider two cases.

i. Assume $k\ne k'$ correspond to different faces $F$ and $F'$. If both are 2-faces, the consistency of the two critical systems with a common $x_2$ and independent coordinates $z,z'$ is equivalent to the consistency of
$$
f|_F=\partial_3f|_F=f|_{F'}=\partial_3f|_{F'}=0.\eqno{(*)}
$$
Indeed, for fixed $x_2$ the monomial map
$$
(x_1,x_3)\longmapsto(z,z')=(x_3x_1^{-k},x_3x_1^{-k'})
$$
is a finite surjective isogeny, because its exponent determinant is $k'-k\ne0$. The four kinetic rows of $(*)$ have rank 4 in every nonvacuous case; if their rank drops, one of the individual face-critical systems is already empty. Thus $(*)$ is an ECI of codimension 4 in three variables and is inconsistent by Theorem \ref{theci}.

If $F$ is a 2-face and $F'$ is a horizontal edge, the same isogeny reduces the question to the three equations
$$
f|_F=\partial_3f|_F=f|_{F'}=0.
$$
If $F'\subset F$ and the kinetic rows have rank 3, their support is two-dimensional, so Theorem \ref{theci} applies; if the rank drops, the individual $F$-critical system is already empty for general $f$. Otherwise the coefficients not shared by the first two equations and the last one make the incidence a proper condition, so it is absent for general $f$. Finally, if $F$ and $F'$ are distinct horizontal edges, their two general one-variable restrictions have no common root.

ii. Assume $k$ and $k'$ correspond to the same horizontal edge $F$, but possibly different edges $E,E'\subset N'_F$. The equation $g_F(x_2)=0$ in the proof of part 2 is the same for every such edge, because it is the restriction of the complete height-$h_F$ layer. Thus $p(Y_k)=p(Y_{k'})$.

The same arguments show that, outside the cases listed in parts 4 and 5, distinct projected local components do not meet along $L$, proving part 6. $\hfill\square$

\subsection{Proof of Lemma \ref{mainl}.4-5: tangent projections.}

\begin{utver}\label{proptangent}
Let $q_i\in\Z_{>0}$, let $z_i=(a,b_i)\in\C^2$, and let
$$
\varphi_{i,j}(x_1,x_2,x_3)=\psi_{i,j}(x_2,x_3)+x_1^{q_i}\eta_{i,j}(x_2,x_3)+x_1^{q_i+1}\rho_{i,j}(x_1,x_2,x_3),
$$
where all the functions involved are holomorphic near the respective points, $\psi_{i,j}(z_i)=0$, $i=1,2$, $j=1,2$. Let $(C_i,z)$ be the image, under the projection forgetting $x_3$, of the germ at $(0,z_i)$ of the curve $\varphi_{i,1}=\varphi_{i,2}=0$, where $z=(0,a)$. Assume that $C_1$ and $C_2$ have no common component. Put
$$
v_i:=\left(\det\left(\begin{smallmatrix}
\partial\psi_{i,1}/\partial x_2(z_i)&\partial\psi_{i,1}/\partial x_3(z_i)\\
\partial\psi_{i,2}/\partial x_2(z_i)&\partial\psi_{i,2}/\partial x_3(z_i)
\end{smallmatrix}\right),
-\det\left(\begin{smallmatrix}
\eta_{i,1}(z_i)&\partial\psi_{i,1}/\partial x_3(z_i)\\
\eta_{i,2}(z_i)&\partial\psi_{i,2}/\partial x_3(z_i)
\end{smallmatrix}\right)\right).
$$
Assume that neither $v_1$ nor $v_2$ is a multiple of $(0,1)$. Then:

1. The intersection number of $C_1$ and $C_2$ at $z$ is at least $\min(q_1,q_2)$.

2. If $q_1=q_2$ and $v_1$ is not a multiple of $v_2$, then the intersection number equals $q_1$.

3. If $q_1<q_2$ and $v_1$ is not a multiple of $(1,0)$, then the intersection number equals $q_1$.
\end{utver}

\begin{proof}
Set $t=x_1$. Since $v_i$ is not a multiple of $(0,1)$, its first coordinate is nonzero, i.e.
$$
\det\left(\begin{smallmatrix}
\partial\psi_{i,1}/\partial x_2(z_i)&\partial\psi_{i,1}/\partial x_3(z_i)\\
\partial\psi_{i,2}/\partial x_2(z_i)&\partial\psi_{i,2}/\partial x_3(z_i)
\end{smallmatrix}\right)\ne0.
$$
Thus, by the implicit function theorem, the germ of the curve $\varphi_{i,1}=\varphi_{i,2}=0$ at $(0,z_i)$ has a unique parametrization
$$
t\longmapsto (t,X_i(t),Y_i(t)),\qquad X_i(0)=a,\quad Y_i(0)=b_i.
$$
In particular, $C_i$ is the graph $x_2=X_i(x_1)$.

Let
$$
A_i:=\left(\begin{matrix}
\partial\psi_{i,1}/\partial x_2(z_i)&\partial\psi_{i,1}/\partial x_3(z_i)\\
\partial\psi_{i,2}/\partial x_2(z_i)&\partial\psi_{i,2}/\partial x_3(z_i)
\end{matrix}\right),\qquad
e_i:=\left(\begin{matrix}\eta_{i,1}(z_i)\\ \eta_{i,2}(z_i)\end{matrix}\right).
$$
The invertibility of $A_i$ and the defining equations first give
$$
X_i(t)-a=O(t^{q_i}),\qquad Y_i(t)-b_i=O(t^{q_i}).
$$
Thus the quadratic Taylor remainder is $O(t^{2q_i})\subset O(t^{q_i+1})$, and substitution into the equations gives
$$
A_i\left(\begin{matrix}X_i(t)-a\\Y_i(t)-b_i\end{matrix}\right)+t^{q_i}e_i=O(t^{q_i+1}).
$$
By Cramer's rule, the first coordinate of this equality is
$$
X_i(t)=a+\frac{(v_i)_2}{(v_i)_1}t^{q_i}+O(t^{q_i+1}).\eqno{(*)}
$$

Since $C_1$ and $C_2$ have no common component, their intersection number at $z$ equals
$$
\operatorname{ord}_t(X_1(t)-X_2(t)).
$$
Formula $(*)$ implies that this order is at least $\min(q_1,q_2)$, proving part 1.

If $q_1=q_2=q$, then
$$
X_1(t)-X_2(t)=\left(\frac{(v_1)_2}{(v_1)_1}-\frac{(v_2)_2}{(v_2)_1}\right)t^q+O(t^{q+1}).
$$
Since the first coordinates of $v_1$ and $v_2$ are nonzero, the coefficient of $t^q$ vanishes if and only if $v_1$ is a multiple of $v_2$. This proves part 2.

Finally, if $q_1<q_2$, then
$$
X_1(t)-X_2(t)=\frac{(v_1)_2}{(v_1)_1}t^{q_1}+O(t^{q_1+1}).
$$
The coefficient of $t^{q_1}$ vanishes if and only if $v_1$ is a multiple of $(1,0)$. This proves part 3.
\end{proof}

We shall also use the following triangular version of this proposition.

\begin{utver}\label{proptriangular}
Let $q_i\in\Z_{>0}$, let $z_i=(a,b_i)\in\C^2$, and let
$$
\varphi_{i,1}(x_1,x_2,x_3)=G(x_1,x_2)+x_1^{q_i}\eta_i(x_2,x_3)+x_1^{q_i+1}\rho_{i,1}(x_1,x_2,x_3),
$$
$$
\varphi_{i,2}(x_1,x_2,x_3)=\frac{\partial\eta_i}{\partial x_3}(x_2,x_3)+x_1\rho_{i,2}(x_1,x_2,x_3),
$$
where all the functions involved are holomorphic near the respective points, $G(0,a)=0$, and, for $i=1,2$,
$$
\frac{\partial G}{\partial x_2}(0,a)\ne0,\qquad
\frac{\partial\eta_i}{\partial x_3}(z_i)=0,\qquad
\frac{\partial^2\eta_i}{\partial x_3^2}(z_i)\ne0.
$$
Let $(C_i,z)$ be the image, under the projection forgetting $x_3$, of the germ at $(0,z_i)$ of the curve $\varphi_{i,1}=\varphi_{i,2}=0$, where $z=(0,a)$. Assume that $C_1$ and $C_2$ have no common component. Put
$$
v_i:=\left(\frac{\partial G}{\partial x_2}(0,a),-\eta_i(z_i)\right).
$$
Then:

1. The intersection number of $C_1$ and $C_2$ at $z$ is at least $\min(q_1,q_2)$.

2. If $q_1=q_2$ and $v_1$ is not a multiple of $v_2$ (equivalently, $\eta_1(z_1)\ne\eta_2(z_2)$), then the intersection number equals $q_1$.

3. If $q_1<q_2$ and $v_1$ is not a multiple of $(1,0)$ (equivalently, $\eta_1(z_1)\ne0$), then the intersection number equals $q_1$.
\end{utver}

\begin{proof}
Set $t=x_1$. At $(0,z_i)$, we have
$$
\det\frac{\partial(\varphi_{i,1},\varphi_{i,2})}{\partial(x_2,x_3)}
=\frac{\partial G}{\partial x_2}(0,a)\frac{\partial^2\eta_i}{\partial x_3^2}(z_i)\ne0.
$$
Thus, by the implicit function theorem, the germ of the curve $\varphi_{i,1}=\varphi_{i,2}=0$ at $(0,z_i)$ has a unique parametrization
$$
t\longmapsto(t,X_i(t),Y_i(t)),\qquad X_i(0)=a,\quad Y_i(0)=b_i.
$$
In particular, $C_i$ is the graph $x_2=X_i(x_1)$.

The implicit function theorem also gives a unique holomorphic function $A(t)$ such that
$$
G(t,A(t))=0,\qquad A(0)=a.
$$
The difference $G(t,x_2)-G(t,A(t))$ equals $x_2-A(t)$ times a holomorphic unit whose value at $(0,a)$ is $\partial G/\partial x_2(0,a)$. Substitution of $(t,X_i(t),Y_i(t))$ into the first equation first gives $X_i(t)-A(t)=O(t^{q_i})$. Also $Y_i(t)-b_i=O(t)$, and hence
$$
\eta_i(X_i(t),Y_i(t))=\eta_i(z_i)+O(t).
$$
It follows that
$$
X_i(t)=A(t)-\frac{\eta_i(z_i)}{\partial G/\partial x_2(0,a)}t^{q_i}+O(t^{q_i+1}).\eqno{(*)}
$$

Since $C_1$ and $C_2$ have no common component, their intersection number at $z$ equals
$$
\operatorname{ord}_t(X_1(t)-X_2(t)).
$$
The common term $A(t)$ cancels in this difference. Formula $(*)$ therefore implies that its order is at least $\min(q_1,q_2)$, proving part 1.

If $q_1=q_2=q$, then
$$
X_1(t)-X_2(t)=-\frac{\eta_1(z_1)-\eta_2(z_2)}{\partial G/\partial x_2(0,a)}t^q+O(t^{q+1}).
$$
Since the first coordinates of $v_1$ and $v_2$ are equal and nonzero, the coefficient of $t^q$ vanishes if and only if $v_1$ is a multiple of $v_2$. This proves part 2.

Finally, if $q_1<q_2$, then
$$
X_1(t)-X_2(t)=-\frac{\eta_1(z_1)}{\partial G/\partial x_2(0,a)}t^{q_1}+O(t^{q_1+1}).
$$
The coefficient of $t^{q_1}$ vanishes if and only if $v_1$ is a multiple of $(1,0)$. This proves part 3.
\end{proof}

We shall apply these propositions to the functions $f_k(x_1,x_2,x_3)$ from the preceding subsection.

\vspace{1ex}

{\it Proof of Lemma \ref{mainl}.4.} The case $p(y)\ne p(y')$ is a tautology. Otherwise, let $y=(y_2,y_3)$ and $y'=(y_2,\varepsilon y_3)$, where $\varepsilon^{d_F}=1$ and $\varepsilon\ne1$. After one monomial normalization, write
$$
\tilde f_k=g(x_2,x_3^{d_F})+x_1^{d/d_F}x_3^p h(x_2,x_3^{d_F})+O(x_1^{d/d_F+1}).
$$
The critical curve is equivalently defined by $(\tilde f_k,\partial_3\tilde f_k)$. Apply Proposition \ref{proptangent} to this same pair at the two points $z_1=y$ and $z_2=y'$. We have $q_1=q_2=d/d_F$, and, up to fixed nonzero factors,
$$
\left(\begin{matrix}v_1\\v_2\end{matrix}\right)=
\left(\begin{matrix}
\partial_2g(y)\,\partial_3^2g(y)&-y_3^ph(y)\,\partial_3^2g(y)\\
\partial_2g(y)\,\partial_3^2g(y)&-\varepsilon^py_3^ph(y)\,\partial_3^2g(y)
\end{matrix}\right).
$$
The systems
$$
g=\partial_3g=\partial_2g=0,\qquad
g=\partial_3g=\partial_3^2g=0,
$$
and
$$
g=\partial_3g=h=0
$$
are full-rank ECI systems with two-dimensional support and hence are inconsistent for general $f$ by Theorem \ref{theci}; in a vacuous thin case the first two equations already have no root. Since $\gcd(p,d_F)=1$, the two vectors are not multiples of each other or $(0,1)$. Their unequal leading coefficients also show that the projected germs have no common component. Proposition \ref{proptangent} gives the intersection number $d/d_F$. $\hfill\square$

\vspace{1ex}

{\it Proof of Lemma \ref{mainl}.5.} The case $p(y)\ne p(y')$ is a tautology, so write $y=(a,b)$ and $y'=(a,b')$. Use the common normalization from the proof of part 2:
$$
\tilde f_k=G_F(x_1,x_2)+x_1^{q_E}H_E(x_2,x_3)+O(x_1^{q_E+1}),\qquad q_E=\frac{d\rho_E}{d_E},
$$
and the analogous formula for $E'$. The function $G_F$ is literally the same in the two formulas. Since $x_3$ is a unit and the critical ideal is saturated by $x_1$, the same curve germ is defined by
$$
\left(\tilde f_k,\,x_3^{-1}x_1^{-q_E}\partial_3\tilde f_k\right).
$$
Its second equation has leading term
$$
x_3^{-1}x_1^{-q_E}\partial_3\tilde f_k
=\frac{\partial H_E}{\partial x_3}+O(x_1).
$$
Thus this pair satisfies Proposition \ref{proptriangular} with $\eta_1=H_E$; the second germ satisfies it with $\eta_2=H_{E'}$. The corresponding vectors are
$$
v_1=\left(\frac{\partial g_F}{\partial x_2}(a),-H_E(a,b)\right),\qquad
v_2=\left(\frac{\partial g_F}{\partial x_2}(a),-H_{E'}(a,b')\right).
$$

For general $f$, the roots of $g_F$ are simple, all critical points of every $H_E(a,\cdot)$ are Morse, and none of their critical values is 0. If $q_E=q_{E'}$, the critical values belonging to distinct germs are also different. For $E=E'$ this follows from Lemma \ref{lmorse}: the coefficients of $g_F$ come from the deleted height layer and are independent of those of $H_E$, so the general map
$$
x_2\longmapsto(g_F(x_2),H_E(x_2,\cdot))
$$
avoids $\{0\}\times M_{B_E}$, which has codimension at least 2. Here
$$
B_E:=\{h(\alpha)-h_F\mid \alpha\in N,\ h(\alpha)\ne h_F,\ \pi_F(\alpha)\in E\}.
$$
Normality gives $\gcd(B_E)=1$ for the edges inherited from $N_F$, while $B_{E_0}=\{-1,1\}$. For $E\ne E'$, equality of two critical values is a proper algebraic condition: after fixing their common coefficients, one may vary a coefficient belonging to one edge and not to the other. Thus it is avoided by general $f$, also when the two edges share an end point.

If $q_E=q_{E'}$, Proposition \ref{proptriangular}.2 now gives intersection number $q_E$. If, after interchanging the two germs, $q_E<q_{E'}$, the nonzero critical value $H_E(a,b)$ allows us to apply Proposition \ref{proptriangular}.3 and gives intersection number $q_E$. In either case this is
$$
\min(q_E,q_{E'})=d\min\left(\frac{\rho_E}{d_E},\frac{\rho_{E'}}{d_{E'}}\right).
$$
The unequal leading coefficients also show that the projected germs have no common component. $\hfill\square$

\begin{lemma}\label{lmorse}
For finite $B\subset\Z$, let $M_B\subset\C^B$ be the closure of the locus of Laurent polynomials having two distinct critical points in $\C^*$ with the same critical value. If $\gcd(B)=1$ (to be distinguished from $\gcd(B-B)=1$ assumed throughout the paper), then $M_B$ has positive codimension. 
\end{lemma}
\begin{proof}
Consider the discriminant $\tilde D$ of the two general polynomials
$$
py+q\quad\mbox{and}\quad ry-\beta(x),\qquad\beta\in\C^B.\eqno{(*)}
$$
By the Cayley trick and $\gcd(B)=1$, this discriminant is reduced, so $\sing\tilde D$ has codimension at least 2. If $\beta$ has two distinct critical points with the same value $c$, then every $(p,q,r)$ satisfying $qr=-pc$ belongs to the component of $\sing\tilde D$ corresponding to two double roots of the eliminated polynomial $p\beta+qr$. This is a two-dimensional family. Hence
$$
\dim M_B+2\le\dim\sing\tilde D\le |B|+1,
$$
and therefore $\dim M_B\le |B|-1$.
\end{proof}

\end{document}